\documentclass[12pt]{amsart}

\usepackage{amsmath}
\numberwithin{equation}{section}
\usepackage{amsthm}
\usepackage{amssymb}
\usepackage{amsfonts}
\usepackage{pinlabel}
\usepackage{graphicx}
\newcommand{\spec}{\mathrm{Spec}}
\newcommand{\torus}{\mathbb{T}}

\newcommand{\ra}{{\rightarrow}}
\newcommand{\eproof}{\hfill\rule{2.2mm}{3.0mm}}

\newcommand{\bvec}{\left [ \begin{array}{c}}
\newcommand{\evec}{\end{array} \right ]}
\newcommand{\Vol}{\mathrm{Vol}}
\newcommand{\R}{{\mathbb R}}

\newcommand{\Z}{{\mathbb Z}}
\newcommand{\T}{{\mathbb T}}
\newcommand{\C}{{\mathbb C}}

\newcommand{\N}{{\mathbb N}}

\newcommand{\dist}{\mathrm{dist}}

\newcommand{\Tu}{\mathbf{T}}

\newcommand{\D}{\noindent}
\newcommand{\Ha}{\mathrm{Vol}}

\newcommand{\grad}{\mathrm{grad}}

\newcommand{\diam}{{\rm diam}}

\renewcommand{\eqref}[1]{(\ref{#1})}
\newcommand{\eat}[1]{}

\newcommand{\bsmall}{\begin{array}[c]{c}}
\newcommand{\esmall}{\end{array}}

\newcommand{\vol}{\mathrm{Vol}}
\newcommand{\Kk}{\mathcal{K}}
\newcommand{\Hh}{\mathcal{H}}
\newcommand{\Eta}{\upsilon}
\newcommand{\Lip}{\mathrm{Lip}}

\renewcommand{\div}{\mathrm{div}}
\theoremstyle{plain}
\newtheorem{theo}{Theorem}[section]
\newtheorem{lem}[theo]{Lemma}
\newtheorem{prop}[theo]{Proposition}

\theoremstyle{definition}
\newtheorem{defi}[theo]{Definition}

\newtheorem{exam}[theo]{Example}
\newtheorem{Claim}[theo]{Claim}

\theoremstyle{remark}
\newtheorem*{rmk}{Remark}

\def\T{\mathbb T}
\def\>{>_{\sigma}}

\title{Sturm-Liouville Estimates for the Spectrum and Cheeger Constant}
\author{Brian Benson}
\address{Department of Mathematics, University of Illinois at Urbana-Champaign, Urbana, IL 61801}
\email{benson9@illinois.edu}

\begin{document}
\baselineskip 18pt

\begin{abstract} Buser's inequality gives an upper bound on the first non-zero eigenvalue of the Laplacian of a closed manifold $M$ in terms of the Cheeger constant $h(M)$. Agol later gave a quantitative improvement of Buser's inequality. Agol's result is less transparent since it is given implicitly by a set of equations, one of which is a differential equation Agol could not solve except when $M$ is three-dimensional. We show that a substitution transforms Agol's differential equation into the Riemann differential equation. Then, we give a proof of Agol's result and also generalize it using Sturm-Liouville theory. Under the same assumptions on $M$, we are able to give upper bounds on the higher eigenvalues of $M$, $\lambda_k(M)$, in terms of the eigenvalues of a Sturm-Liouville problem which depends on $h(M)$. We then compare the Weyl asymptotic of $\lambda_k(M)$ given by the works of Cheng, Gromov, and B\'erard-Besson-Gallot to the asymptotics of our Sturm-Liouville problems given by Atkinson-Mingarelli.
\end{abstract}

\maketitle

\section{Introduction} \label{section:intro}
\subsection{Summary of Results}

We give an upper bound on the eigenvalues of the Laplacian on a compact Riemannian manifold in terms of the Cheeger constant of the manifold, denoted $h(M)$. Buser was the first to give such an inequality for the first non-zero eigenvalue of the manifold, denoted $\lambda_1(M)$ \cite{B82}. Agol recently gave a quantitative improvement of Buser's inequality \cite{IA}. The drawback of Agol's improvement is that it is given implicitly by a set of equations, one of which is a differential equation that Agol could only solve in the case of 3-manifolds. We show that a substitution transforms Agol's differential equation into the Riemann differential equation, which is well understood.

We use Sturm-Liouville theory as a framework for giving upper bounds on the spectrum of the manifold $M$ in terms of $h(M)$. This allows us to not only replicate the known bounds on $\lambda_1(M)$ in terms of $h(M)$, but to extend these results to give upper bounds for the higher eigenvalues, denoted $\lambda_k(M)$, in terms of $h(M)$. To our knowledge, these are the first upper bounds for $\lambda_k(M)$ in terms of $h(M)$. Our bounds are eigenvalues of one-dimensional Sturm-Liouville problems which depend on the parameter $h=h(M)$. A consequence of Sturm-Liouville theory is that these bounds are differentiable almost everywhere as functions of $h>0$. We also consider asymptotic growth rates for these upper bounds in terms of $k$ and compare them to known asymptotic growth rates for $\lambda_k(M)$.

\begin{figure}[htb]
\labellist
\small\hair 2pt
 \pinlabel {$\lambda_1(h)$} [ ] at 2 226
 \pinlabel {$h$} [ ] at 370 13
 \pinlabel {Agol} [ ] at 239 66
 \pinlabel {Buser} [ ] at 185 123
\endlabellist
\centering
\includegraphics[scale=1.0]{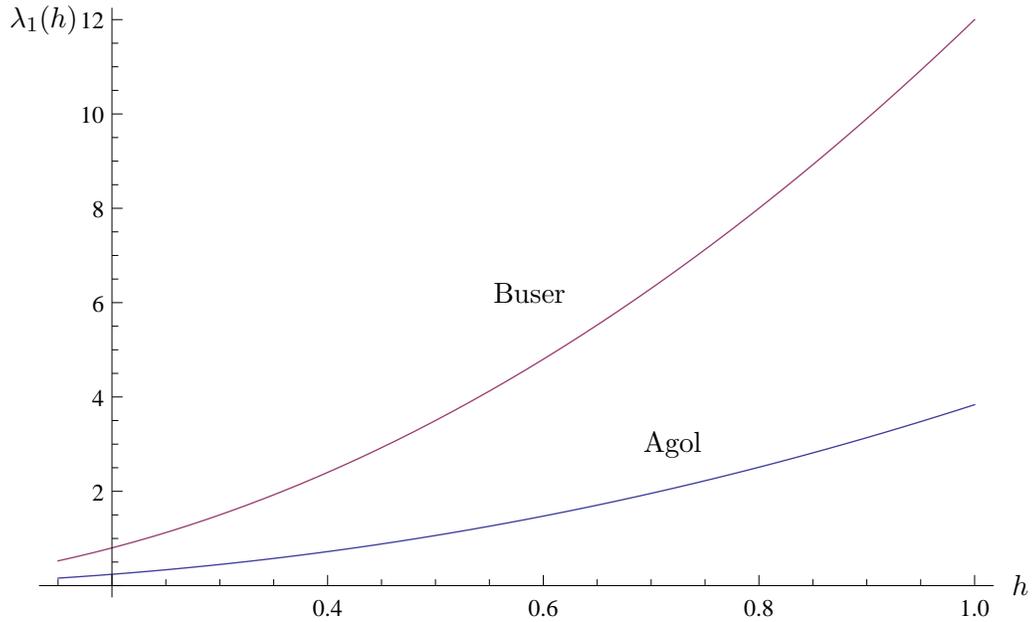}
\caption{Buser's inequality \cite{B82} versus Agol's improvement \cite{IA} in dimension 2 for $\lambda_1$ as a function of $h$. In both estimates, Ricci curvature is bounded from below by $-(n-1)$.}
\label{fig:PL}
\end{figure}

\begin{figure}[htb]
\labellist
\small\hair 2pt
 \pinlabel {$h$} [ ] at 370 13
 \pinlabel {$\lambda_1(h)$} [ ] at 13 226
 \pinlabel {$n=2$} [ ] at 370 193
 \pinlabel {$n=3$} [ ] at 370 203
 \pinlabel {$n=4$} [ ] at 370 213
 \pinlabel {$n=5$} [ ] at 370 223
 \pinlabel {$n=6$} [ ] at 370 233
\endlabellist
\centering
\includegraphics[scale=1.0]{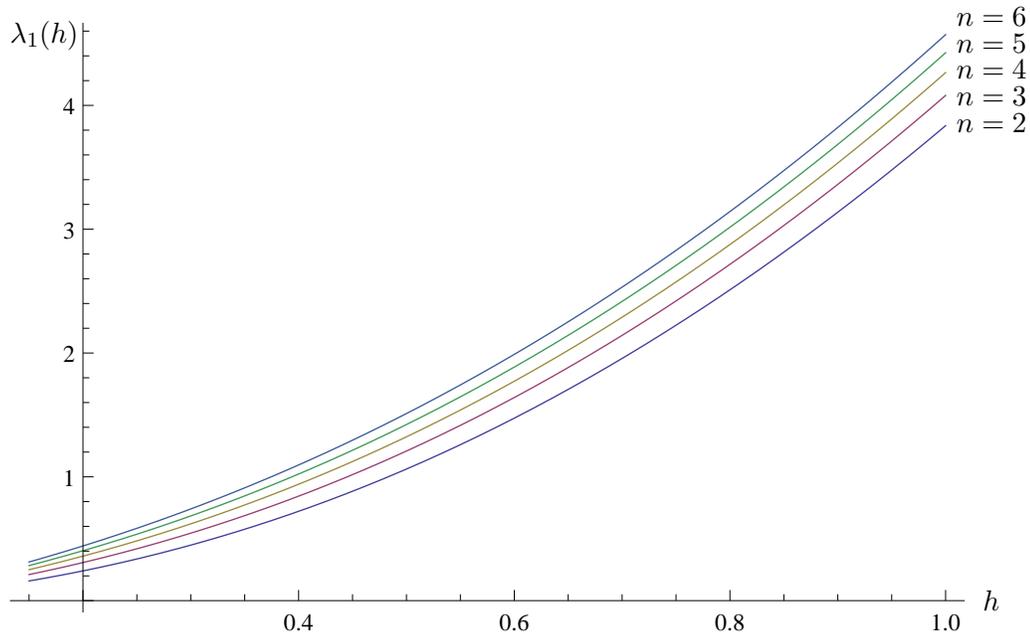}
\caption{Agol's improvement for $\lambda_1$ as a function of $h$ for dimensions $n=2,3,4,5,6$.  }
\label{fig:nchange}
\end{figure}

For additional motivation, here are two plots of numerical data corresponding to Agol's improvement. Specifically, Figure \ref{fig:PL} gives a comparison between Buser's inequality and Agol's improvement in dimension 2. Figure \ref{fig:nchange} shows Agol's upper bound on $\lambda_1$ as a function of $h$ for dimensions $n=2,3,4,5,6$ to demonstrate that plots for higher dimensions are similar up to scale.

\subsection{Notation and Conventions} Let $M$ be a closed $n$-dimensional Riemannian manifold, with $n\geq 2$.  For $u \in C^2(M)$, the geometer's Laplacian of $u$ is $\Delta u = - \div \big (\grad (u)\big )$.  Eigenvalues of the Laplacian are real values $\lambda$ such that $\Delta u=\lambda u$ for some $u \in C^2(M)$ where $u$ satisfies the Dirichlet boundary condition $u|_{\partial M}=0$. 
For closed manifolds, the spectrum starts with $\lambda_0$:
   $$0=\lambda_0 < \lambda_1 \leq \lambda_2 \leq \lambda_3 \leq \cdots,$$
while for manifolds with boundary, the spectrum begins with $\lambda_1$:
   $$0 < \lambda_1 \leq \lambda_2 \leq \lambda_3 \cdots.$$
In both situations, the first positive eigenvalue is denoted $\lambda_1$.
We will study the connections between the spectrum of the Laplacian of a manifold $M$ and its Cheeger constant, defined as follows.
\begin{defi}\label{defi:Cheeger}
The Cheeger constant of a closed $n$-dimensional Riemannian manifold $M$ is defined as 
		 $$h(M) = \inf \frac{\Vol_{n-1}(\partial A)}{\Vol_n (A)}$$ where $A \subset M$ and $\partial A$ is a smooth codimension-1 submanifold of $M$ and
		 $\Vol_n(A) \leq \frac{1}{2}\Vol_n(M)$.  The quantity $\Vol_{n-1} (\partial A) / \Vol_n (A)$ is called the isoperimetric ratio of the set $A$.
\end{defi}

\subsection{Historical Motivation}  Cheeger proved that $\lambda_1(M) \geq \left (h(M)/2 \right )^2$, providing the initial motivation for defining the Cheeger constant \cite{C69}.  However, even before the work of Cheeger, the classical isoperimetric inequality gave the following result for subsets  of the $n$-sphere. For any $U \subset S^n$, choose an open ball $O \subset S^n$ so that $\Vol_n(O)=\Vol_n(U)$.  Then the classical isoperimetric inequality can be stated as
	$$\Vol_{n-1} (\partial U) \geq \Vol_{n-1} (\partial O).$$
L\'evy \cite{L51} extended the classical isoperimetric inequality to the case of convex hypersurfaces in $\R^{n+1}$. Later Gromov showed that L\'evy's method can be canonically extended to the case of closed Riemannian manifolds \cite[Appendix]{Grom2}. In particular, Gromov proved that when the Ricci curvature of $M$ is bounded below by $-(n-1)$ and $d=\diam (M)$, then  $\lambda_1(M) \geq e^{-2(n-1)d}$ using L\'evy's method and Cheeger's inequality from above; this result was also proved independently by Li and Yau \cite{LY80}. In addition, for any $\epsilon >0$, letting $N=N(\epsilon)$ be the minimum integer such that $M$ can be covered by $N$ balls of radius $\epsilon$, Gromov showed that there exist positive constants $C_1, C_2$ depending on the lower bound on Ricci curvature such that 
	$$\epsilon^{-2}C_1^{1+\epsilon} \leq \lambda_N (M) \leq \epsilon^{-2} C_2^{1+\epsilon}.$$  In summary, Gromov was able to obtain bounds on higher eigenvalues of $M$ by taking $\epsilon$ to be small.

Interestingly, Kr\"oger gave lower bounds on $\lambda_1(M)$ in terms of the eigenvalues of a corresponding Sturm-Liouville problem depending on the dimension, Ricci curvature, and diameter of $M$. In addition, he gave examples where his estimates are sharper than estimates given by Cheeger's inequality \cite{K92}\cite{K97}.

Buser, citing Gromov's work as motivation, proved that for $M$ a closed Riemannian manifold with Ricci curvature bounded below by $ -\delta^2(n-1)$, then $$\lambda_1 (M) \leq 2\delta (n-1) h(M)+10h^2(M)$$ \cite{B82}. Combining the results of Buser and Cheeger, we have the following qualitative statement: For closed manifolds, the first eigenvalue of the Laplacian is controlled by the Cheeger constant.  

Agol observed that Buser's inequality gave a far from sharp estimate for certain hyperbolic 3-manifolds, motivating him to improve it \cite{IA}. Agol uses a function $\bar{J}$ and parameter $\bar{T}$ depending on the dimension $n$ and Cheeger constant $h(M)$ as follows. The function $\bar{J}$ is given by 
\begin{equation} \label{eq:barJ} \bar{J}(\tau)=\left  (\cosh(\tau)+\frac{h}{n-1}\sinh(\tau) \right )^{\! \! \! \frac{n-1}{2}}.\end{equation}
Further, $\bar{T} \in (0, \infty)$ is defined implicitly by the equation $$\frac{1}{h} = \int_0^{\bar{T}} \bar{J}^2(\tau) \, d \tau ,$$ which is valid for every $h$ since the right hand side can take on all values from $0$ to $\infty$ because the integral approaches $0$ as $\bar{T} \to 0$ and approaches $\infty$ as $\bar{T} \to \infty$. Agol proved the following:

\begin{theo}\label{theo:Agol} {\bf (Agol \cite{IA})} There is a function $\lambda (h)$ such that for all closed Riemannian $n$-manifolds $M$ with Ricci curvature bounded from below by $-(n-1)$ we have that $\lambda\big (h(M)\big ) \geq \lambda_1(M)$. Moreover, we can take $\lambda (h)$ to be the least positive number such that there exists a (non-trivial) $y \in C^{\infty}[0,\bar{T}]$ satisfying
			\begin{equation*} %\label{eq:AgolSL} 
				y''=\left (\frac{\bar{J}''}{\bar{J}} - \lambda (h) \right ) y, \qquad y(0)=0, \qquad 									y'\left (\bar{T}\right )=\bar{J}'\left (\bar{T}\right ), \qquad y\left (\bar{T}\right )=\bar{J}\left (\bar{T}\right ).
			\end{equation*}
\end{theo}

Far less is known about the relationship between $\lambda_k(M)$ and $h(M)$ for a closed manifold. However, asymptotic results for the higher eigenvalues of $M$ in terms of the dimension, Ricci curvature, and volume of the manifold have been thoroughly developed. Specifically,  B\'erard, Besson, and Gallot \cite{BBG85}, building on Cheng \cite{C75a, C75b} and Gromov \cite[Appendix]{Grom2}, showed that, for $M$ closed with Ricci curvature bounded from below by $R$, the asymptotic of $\lambda_k (M)$ is of the order $k^{2/n}$ where the constant depends on $n$, $R$, $\Vol_n(M)$.

\subsection{Detailed Description of Results}\label{sec:detail}
In section \ref{sec:Riemann}, we show that a substitution transforms Agol's differential equation into the Riemann differential equation. In particular, we give a restatement of Agol's Thorem \ref{theo:Agol} as our Theorem \ref{theo:SLE}.

Let $\Eta (h) >0$ be an upper bound on the absolute value of the mean curvature of $\Sigma^0$, the smooth part of a rectifiable current $\Sigma$ dividing $M$ into two sets $A$ and $B$ with $A \cup B =M$ and $A \cap B = \Sigma$ so that $h(M)=\Vol_{n-1}(\Sigma)/\min \{ \Vol_n(A), \Vol_n(B) \}$. Let $q=q(n,\lambda)=\frac{n-1-2\lambda}{2}, r=r(n)=\frac{(n-1)(n-3)}{4},$ and $s=q+r$.  Further, let $a=\frac{\Eta +1}{\Eta -1}$, and $b=e^{2\bar{T}}\frac{\Eta +1}{\Eta -1}$. Under these assumptions, Agol's Theorem \ref{theo:Agol} is equivalent to the following theorem by taking $\Eta(h)=\frac{h}{n-1}$ when $h \neq n-1$.\footnote{Note that when $h=n-1$, Agol's differential equation simplifies greatly using the identity $\cosh(\tau)+\sinh (\tau)=e^{\tau}$.}  The proof of this Theorem appears in section \ref{sec:Riemann}.
\begin{theo}\label{theo:Main3} There exists a function $\lambda(h)$ such that for all closed Riemannian $n$-manifolds $M$ with Ricci curvature lower bound of $-(n-1)$ we have that $\lambda (h(M)) \geq \lambda_1(M)$. Moreover, we can take $\lambda(h)$ to be the smallest positive number $\lambda$ such that there exists $y:[a,b] \ra \R$ (or $y:[b,a]\ra \R$ when $b<a$) satisfying
			\begin{equation}\label{eq:Riem2} y''(z) + \frac{1}{z} y'(z) - 
				\frac{q(n, \lambda)(z-1)^2+r(n)(z+1)^2}{4z^2(z-1)^2}y(z)=0,
			\end{equation} 
			\begin{equation*}
				y(a)=0,\qquad y'(b)=\frac{\bar{J}'(T)}{2b},\qquad y(b)=\bar{J}(T).
			\end{equation*}
		Finally, (\ref{eq:Riem2}) is an example of the Riemann differential equation with regular singularities at $0,1, \infty$ and respective local exponents 
		$$\begin{array}{c}\left \{ -\frac{\sqrt{q+r}}{2}, \frac{\sqrt{q+r}}{2} \right \}, \left \{\frac{1-\sqrt{1+4r}}{2},\frac{1+\sqrt{1+4r}}{2} \right \}, \left \{ -\frac{\sqrt{q+r}}{2}, \frac{\sqrt{q+r}}{2} \right \}.\end{array}$$  
Therefore, the solutions of (\ref{eq:Riem2}) for any $z \in \C$ are given by branches of the Riemann $P$-function.
\end{theo}

\begin{rmk}
 Since $n\geq 2$, it follows that $1+4r \geq 0$; however, it is possible that $q+r<0$, and so $\sqrt{q+r}$ has a non-zero imaginary component. In all but the simplest of cases, one should think of the variable $z$ in equation (\ref{eq:Riem2}) as lying in the complex plane. Then the function $y:[a,b]\ra \R$ (or $y:[b,a]\ra \R$) is the real part of a branch of the multi-valued function $y(z)$ given by equation (\ref{eq:Riem2}) when $z \in \C$.\footnote{These branches are hypergeometric functions which we did not find to be very practical in giving numerical upper bounds for $\lambda_1(M)$ in terms of $h(M)$. This is one reason for adopting the point of view of Sturm-Liouville theory.}
\end{rmk}

We then consider the approaches of Buser \cite{B82} and Agol \cite{IA} within the framework of Sturm-Liouville theory.  In section \ref{sec:AgolRiem}, we provide a proof of Agol's Theorem \ref{theo:Agol} which uses the spectral theorem in place of the variational principle used by Agol \cite{IA}.  This new approach allows us to give upper bounds on higher eigenvalues of $\Delta$ in terms of $h(M)$ by using a Sturm-Liouville problem. Like in Buser and Agol, we assume that $M$ is closed with Ricci curvature bounded below by $-(n-1)\delta^2$ for $\delta \geq 0$. 

We use the notation
	$$s_{\delta}(\tau) := \left \{ \begin{array}{rr} \frac{\sinh (\delta \tau)}{\delta},& \delta >0,\\ \tau, & \delta=0,\end{array}\right. $$
	$$c_{\delta}(\tau) :=\frac{ds_{\delta}}{d\tau}. $$
  For any real number $t \geq 0$ and $\delta \geq 0$, define 
	\begin{equation}\label{eq:Jdel}
		J_{\delta} (\tau,t) := \big (c_{\delta}(\tau) + t s_{\delta}(\tau) \big )^{n-1}. 
	\end{equation}
Define weight functions $p$ and $w_i$ for $i=1,2$ which depend on $h$, by
	\begin{gather*}
		p(\tau)=w_1(\tau)=J_{\delta}(\tau, \Eta)= \big ( c_{\delta}(\tau) + \Eta s_{\delta} (\tau) \big )^{n-1}, \\
		w_2(\tau)=1 - h\int_0^{\tau} J_{\delta} \big (x, \Eta (h) \big ) \, dx.
	\end{gather*}
We define $T$ implicitly by
	\begin{equation*}
		\frac{1}{h} = \int_0^T J_{\delta} \big (\tau, \Eta (h) \big ) \, d\tau.
	\end{equation*}
As with $\bar{T}$, the implicit definition of $T$ is valid for any $h$ because the integral approaches $0$ as $T \to 0$ and approaches $\infty$ as $T \to \infty$. Also, the weight functions are all positive on the closed interval $[0,T]$, except for $w_2$ which degenerates to $0$ at $\tau =T$.
	For $i=1,2$, we consider the formally self-adjoint differential operator $L_i$ given by $$L_i u=-\frac{1}{w_i(\tau)}\frac{d}{d\tau}\left ( p(\tau ) \frac{du}{d\tau}\right )$$
and let $\xi(k) := \left \lceil \frac{k+1}{2}\right \rceil$.
For $h=h(M)$, let $\omega_i(h)$ be the regular Sturm-Liouville problem given by 
		\begin{equation}\label{eq:SLPu}
			L_i u =\lambda u, \qquad u(0)=0, \qquad u'(T)=0,
		\end{equation}
for a function $u$ in a suitable Sobolev space to be defined in Section \ref{sec:SLP}. Denote the $k$-th eigenvalue of $\omega_i(h)$ by $\lambda_k\big (\omega_i(h)\big )$.
 In section \ref{sec:SLP}, we prove the following generalization of Agol's Theorem \ref{theo:Agol}:

\begin{theo}\label{theo:SLE}
	Let $T$, $M$, and $\omega_i(h)$ be as above with $h=h(M)$. Then 
\begin{equation}\label{eq:SLE}
\lambda_1(M) \leq \lambda_1\big (\omega_1(h)\big ), \qquad \lambda_k(M) \leq \lambda_{\xi(k)} \big (\omega_2(h)\big ).
\end{equation} 
\end{theo}
\begin{rmk}
The Sturm-Liouville problem $\omega_2$ does not give as sharp of a bound for $\lambda_1$ compared to the Sturm-Liouville problem $\omega_1$; in other words, $\lambda_1\big (\omega_1(h)\big ) \leq \lambda_1\big (\omega_2(h)\big )$ for each $h \in (0, \infty)$.
\end{rmk}

\begin{figure}[htb]
\labellist
\small\hair 2pt
 \pinlabel {$h$} [ ] at 365 4
 \pinlabel {$\lambda_1\big (\omega_2(h)\big )$} [ ] at 337 28
 \pinlabel {$\lambda_2\big (\omega_2(h)\big )$} [ ] at 325 63
 \pinlabel {$\lambda_3\big (\omega_2(h)\big )$} [ ] at 318 130
 \pinlabel {$\lambda_4\big (\omega_2(h)\big )$} [ ] at 290 200
\endlabellist
\centering
\includegraphics[scale=1.0]{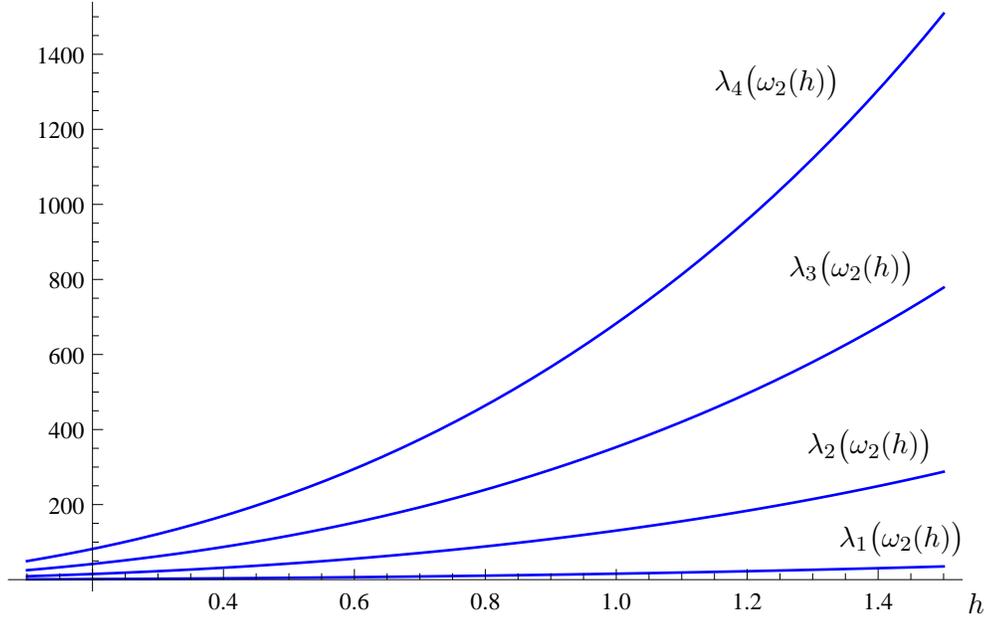}
\caption{The functions $\lambda_k\big (\omega_2(h)\big )$ given by the Sturm-Liouville problem $\omega_2(h)$ for $\delta=1$, $n=2$, and $k=1,2,3,4$.}
\label{fig:HigherEigen}
\end{figure}

See Figure \ref{fig:HigherEigen} for an example of the bounds on higher eigenvalues given by $\omega_2(h)$. By Theorem \ref{theo:SLE}: $\lambda_1(M) \leq \lambda_1\big (\omega_2(h)\big )$, $\lambda_2(M)\leq \lambda_3(M) \leq \lambda_2\big (\omega_2(h)\big )$,  $\lambda_4(M) \leq \lambda_5(M) \leq \lambda_3\big (\omega_2(h)\big )$, and $\lambda_6(M)\leq \lambda_7 (M) \leq \lambda_4\big (\omega_2(h)\big )$.

While one might guess that the Sturm-Liouville problem arising from Agol's work can be extended to higher eigenvalues of $M$ as a direct consequence of our use of the spectral theorem, this is not the case. The proof that $\lambda_1(M) \leq \lambda_1\big ( \omega_1(h) \big )$, which is equivalent to Agol's Theorem \ref{theo:Agol}, uses the fact that the eigenfunction corresponding to the eigenvalue of $\omega_1(h)$ is monotone. While this is certainly true for the eigenfunction corresponding to $\lambda_1$, linear independence of eigenfunctions in $L^2$ means that this cannot hold for eigenfunctions corresponding to larger eigenvalues. Thus, we cannot extend $\omega_1(h)$ directly to give upper bounds for $\lambda_k (M)$ when $k >1$. This is the reason we must use a second Sturm-Liouville problem to give upper bounds for $\lambda_k(M)$ in terms of $h(M)$ for $k\geq 1$.

From the works of Cheng \cite{C75a, C75b}, Gromov \cite[Appendix]{Grom2}, and B\'erard, Besson, and Gallot \cite{BBG85}, it is known that the Weyl law $\lambda_k(M) \asymp k^{2/n}$ holds. A consequence of the Sturm-Liouville framework is that we can apply the work of Atkinson and Mingarelli \cite{AM87}, which gives the following as an immediate consequence:
		There exists a constant $\widetilde{C}=\widetilde{C}(h,n, \delta)$ so that 
			$$\frac{\lambda_{j}\big (\omega_2(h)\big )}{j^2} \, \longrightarrow \, \widetilde{C}$$ as $j \ra \infty$.
		Specifically, we take $\widetilde{C}=\pi^2\left ( \int_0^{T}	\sqrt{\frac{w_2}{p}}\,d\tau \right )^{\!\!-2}$.
	Since $\xi(k)^2$ grows like $k^2$, which is faster than $k^{2/n}$, this approach does not give sharp quantitative upper bounds on $\lambda_k (M)$ for large $k$. In fact, the following example shows that, when the Cheeger constant of $M$ replaces the $n$-volume of $M$ as input data, an upper bound cannot satisfy the Weyl law in general.
	
\begin{exam}\label{exam:flattori}
	We consider examples of flat tori $\T_i$ which are given by the quotient $\R^2 / A_i \Z^2$ for $i \in \N$ where $A_i = \left ( \begin{array}{cc} \frac{1}{2^{i-1}} & 0\\ 0 & 1 \end{array} \right )$. It follows that the Cheeger constant $h(\T_i) = 4$ for each $i$. The eigenvalues of $\T_i$ are of the form $\lambda(\T_i)=4\pi^2 \left ( (2^{i-1}x_1)^2 + x_2^2 \right )$ for $x_1, x_2 \in \Z$.
	
	We would like to see how $\lambda(\T_i)$, when ordered and indexed by $k \in \N$, grows compared to $C\cdot k$ for any fixed constant $C>0$. To see this, we will use a lattice counting argument. Let $L$ be the lattice generated by $A_i^{-1}\Z^2$. Define $p(x,y)=x^2+y^2$. Order the points in $L$ using $p$. Then our question becomes: What is $p$ of the $k$-th point of $L$ using this ordering? A second way of viewing this is that given $P \in \N$, we want to know the number of points $(x,y) \in L$ with $p(x,y) <P$. The answer is the number of lattice  points in the circle of radius $\sqrt{P}$, which is $$\frac{\Vol_2\left (D_{\sqrt{P}}^2\right )}{\text{Covolume of }L}=\frac{\pi P}{2^{i-1}},$$ since the covolume of $L$ is $2^{i-1}$. So the $k$-th point is mapped to approximately $2^{i-1}k/\pi$ under the image of $p$. So the $k$-th eigenvalue of $\T_i$ is asymptotic to $4\pi(\text{Covolume of }L)k = 2^{i+1}k$. Now, $2^{i+1} \nearrow \infty$ as $i \ra \infty$. So we conclude that, in general, there does not exist a constant $C=C\big (n,\delta, h(M)\big )$ so that $\lambda_k(M) \leq Ck^{2/n}$.
\end{exam}

Interestingly, Miclo recently gave a generalized Cheeger constant for each $k \in \N$ which generalizes both Cheeger's and Buser's inequalities corresponding to $\lambda_k(M)$ \cite{LM}. While there is a different constant for each eigenvalue, these bounds satisfy the Weyl law above.  

\subsection{Outline of Proof of Theorem \ref{theo:SLE}}\label{sec:outline}

\begin{figure}[htb]
\labellist
\small\hair 4pt
 \pinlabel {$\Sigma$} [ ] at 205 100
 \pinlabel {$d^{-1}(\tau)$} [ ] at 335 80
  \pinlabel {$D$} [ ] at 390 195
\endlabellist
\centering
\includegraphics[scale=0.7]{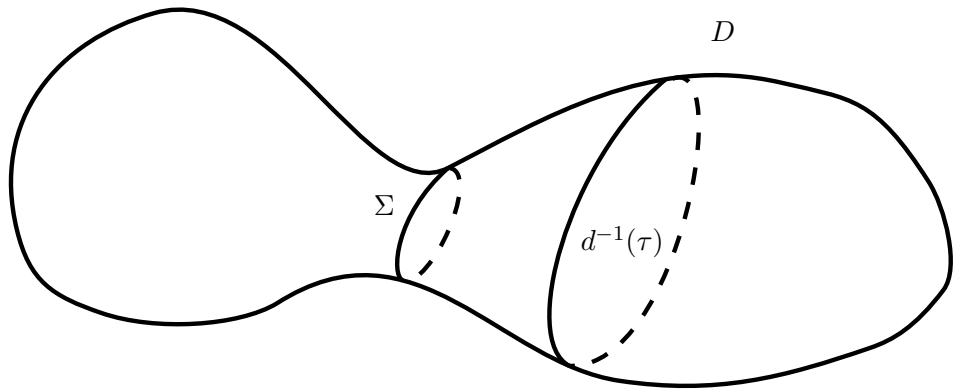}
\caption{Visual example of the level sets $d^{-1}(\tau)$ of $\Sigma$ on $M$ for $\tau>0$.}
\label{fig:Level}
\end{figure}

For convenience, here is a short outline of the proof of Theorem \ref{theo:SLE}:
\begin{enumerate}
\item There is a rectifiable current $\Sigma$ of dimension $n-1$ whose isoperimetric ratio realizes the Cheeger constant.
\item Fix $k \in \N$. Take $D$ to be the closure of the component of $M-\Sigma$ where $\lambda_k(D) \geq \lambda_k(M-D)$ with respect to the Dirichlet problem (eigenfunctions vanish on $\Sigma$).
\item Estimate $\lambda_k(M)$ from above using Proposition \ref{prop:Courant} with $l=k$, which says that $$\lambda_{2k-1}(M) \leq \lambda_k(D).$$
\item The Poincar\'e minimax principle gives the following Rayleigh quotient: 
	\begin{equation}\label{eq:Poin} \lambda_k (D) = \inf_V \sup_{f \in V,\, f \not\equiv 0} \frac{\int_D \| \grad f \|^2 \, d\Vol_n}{\int_D f^2 \, d\Vol_n} \end{equation}
where $V$ runs over all $k$-dimensional subspaces of $H_0^1(D)$, when $\partial D \neq \emptyset$.
\item Take $d:M \ra \R$ to be the signed distance to $\Sigma$, where $d^{-1}[0,\infty)=D$.
\item We estimate (\ref{eq:Poin}) from above by a one-dimensional test function defined which is constant on the level sets $d^{-1}(\tau)$ off of $\Sigma$:
	\begin{equation}\label{eq:Ray}
		\lambda_k(D) \leq \inf_{V_{\ast}} \sup_{g \in V_{\ast}} 
			\frac{\int_0^{\infty}\big (g'(\tau)\big )^{\! 2} \Vol_{n-1} \big (d^{-1}(\tau)\big )d\tau}{\int_0^{\infty} g^2(\tau) \Vol_{n-1} \big (d^{-1}(\tau)\big ) d\tau},
	\end{equation}
where $V_{\ast}$ ranges over $k$-dimensional subspaces of $H^1[0, \infty)$. See Figure \ref{fig:Level} for a visual example of $\Sigma$ and $d^{-1}(\tau)$.
\item Heintze and Karcher \cite{HK78} give a scaling function $J(\tau)$ such that  
\begin{equation}\label{eq:UBlevel}\Vol_{n-1}\big (d^{-1}(\tau)\big ) \leq \Vol_{n-1}(\Sigma)J (\tau),\end{equation}
for $\tau \in [0,\infty)$; see Lemma \ref{lem:HK2}.
\item There exists $T>0$ so that restricting (\ref{eq:Ray}) to the class of test functions such that $g (\tau)=g(T)$ for all $\tau >T$, and combining with (\ref{eq:UBlevel}), we have 
\begin{equation}\label{eq:Ray2}
\lambda_1(D) \leq \frac{\int_0^{T} \big (g'(\tau)\big )^{\! 2} J(\tau) \, d\tau}{\int_0^{T} g^2(\tau) J(\tau) \, d\tau}.
\end{equation}
 This uses the fact that $g$ can be taken to be monotone in the minimization of the quotient on the right hand side of (\ref{eq:Ray}).
\item For higher eigenvalues of $\lambda_k(D)$, we provide a lower bound for $\Vol_{n-1}\big (d^{-1}(\tau)\big )$ in Lemma \ref{lem:LSigTau}, namely: 
\begin{equation}\label{eq:LBlevel}
\Vol_{n-1}\big (d^{-1}(\tau)\big ) \geq \Vol_{n-1}(\Sigma)\left ( \frac{\Vol_n \big ( d^{-1}(\tau, \infty) \big )}{\Vol_n (D)}\right )\end{equation}
for $\tau \in [0,T]$ almost everywhere.
\item Applying the upper bound (\ref{eq:UBlevel}) and lower bound (\ref{eq:LBlevel}) for $\Vol_n\big (d^{-1}(0,\tau)\big )$ to the right side of  (\ref{eq:Ray}), gives 
\begin{equation}\label{eq:Ray3}
\lambda_k(D) \leq \inf_{V_{\ast}} \sup_{g \in V_{\ast}} \frac{\int^{T}_0 \big (g'(\tau)\big )^{\! 2} J (\tau) \, d\tau}{\int_0^{T} g^2(\tau) \left (1 - h\int_0^{\tau} J (x) \, dx \right ) \, d\tau}.
\end{equation}
where $V_{\ast}$ ranges over $k$-dimensional subspaces of $H^1[0, T]$.
\item The spectral theorem is applied to show that the test functions $g$ on the right sides of (\ref{eq:Ray2}) and (\ref{eq:Ray3}) are exactly the solutions $u$ of the respective Sturm-Liouville problems in (\ref{eq:SLPu}); see Lemma \ref{lem:Main2}. The result follows.
\end{enumerate}

\pagebreak
\subsection{Plots}
Bailey, Everitt, and Zettl give a Fortran program called {\em SLEIGN2} which estimates the eigenvalues of Sturm-Liouville problems \cite{BEZ}. We use this program in coordination with Mathematica to produce the plots of $\lambda_k \big (\omega_i(h)\big )$ seen herein.

\subsection{Acknowledgments} The author acknowledges support from National Science Foundation grant DMS 0838434 ``EMSW21MCTP: Research Experience for Graduate Students''.  The author would like to thank Ian Agol for his mentorship, allowing him access to his unpublished work, and for permitting him to reference this work freely herein.  The author thanks his advisor Nathan Dunfield for helping produce the plots, suggesting the works of Agol and Buser, and for invaluable suggestions and insights.  The author also thanks Pierre Albin, Bruno Colbois, Neal Coleman, Chris Judge, Gabriele La Nave, Richard Laugesen, Jeremy Tyson, and Brian White for helpful discussions and emails. The author is grateful to Richard Laugesen for providing Example \ref{exam:tent}. In addition, the author thanks Marwa Balti for helpful comments on the writing of this manuscript.

%%%%%%%%%%%%%%%%%%%%%%%%%%%%%%%%%%%%%%%%%%%%%%%%%%%%%%%%%%%%%%%%%%%%%%%%%%%%%%%%%%%%%%%%%%%%%%%%%%%%%%%%%%%%%%%%%%%%%%%%%%%%%%%%%%%
\section{Eigenvalues, the Cheeger Constant, and Minimizing Currents}
\label{sec:eigen}
In this section, we bound $\lambda_k(M)$ from above by a Rayleigh quotient which uses a test function that is constant on each level set $d^{-1}(\tau)$. To give bounds on $\lambda_1(M)$, it suffices to bound $\Vol_{n-1}\big (d^{-1}(\tau)\big )$ from above by $\Vol_{n-1}(\Sigma)$ multiplied by a smooth scaling function which depends on $h(M)$; such a result follows from the work of Heintze and Karcher \cite{HK78}. To give bounds on the higher eigenvalues, $\lambda_k(M)$, we will also bound $\Vol_{n-1}\big (d^{-1}(\tau)\big )$ from below by $\Vol_{n-1}(\Sigma)$ multiplied by a scaling function depending on $h(M)$.

\subsection{Separating Rectifiable Currents} 
Buser, using Almgren's work \cite{FA76}, showed that whenever $M$ is closed, there is a closed set $A \subseteq M$ with $\Vol_n (A) \leq \frac{1}{2}\Vol_n (M)$ such that the isoperimetic ratio of $A$ realizes $h(M)$ \cite{B82}.  Moreover, $\Sigma=\partial A$ is a rectifiable current of codimension-1 in $M$, see Section \ref{sec:current} for a definition. For dimensions $n\leq 7$, Morgan \cite{FM03} showed that $\Sigma$ is a smooth submanifold. For an overview of why $\Sigma$ need not be a hypersurface for dimensions $n=8$ and higher, see Federer \cite{HF69} and Morgan \cite{FM09}.

The fact that $\Sigma$ may not be a smooth hypersurface will not cause too much concern. As Gromov points out, with the help of Almgren's work \cite{FA76}, if $x \in M$ and $\gamma$ is a geodesic segment from $x$ to $\Sigma$ realizing $\dist(x,\Sigma)$, then $\gamma$ ends at a nonsingular point of $\Sigma$ \cite[Appendix]{Grom2}. Building on the work of Federer \cite{HF70} and Almgren \cite{FA76}, Morgan proved that $\Sigma$ is locally a smooth $C^{\infty}$-submanifold of $M$ except for a set of Hausdorff dimension at most $n-7$ \cite{FM03}. Thus $\Sigma^0$ is a smooth hypersurface and $\Vol_{n-1}(\Sigma)=\Vol_{n-1}(\Sigma^0)$. Finally, it is well-known that $\Sigma^0$ must have constant mean curvature; see Ros \cite{R05}.

We will divide $M$ into two sets $A$ and $B$, with $\vol_n(A) \leq \vol_n(B)$, via a rectifiable current $\Sigma$ so that $A \cup B = M$ and $A \cap B = \Sigma$ and so that 
\begin{equation*}%\label{eq:hAB}
h(M) = \frac{\vol_{n-1} (\Sigma)}{\vol_n (A)}.
\end{equation*}

\subsection{Minimax Principles}
We now show how to give an upper bound for $\lambda_k (M)$ in terms of a Rayleigh quotient on a space of functions defined on a compact interval of the real line. The methods we use closely follow the arguments given in Buser \cite{B82} for $\lambda_1(M)$. We begin to generalize to $\lambda_k (M)$ by applying the Poincar\'e minimax principle.

We use the decomposition of $M$ into the components $A$ and $B$ to give an upper bound on the eigenvalues of $M$ in terms of Dirichlet eigenvalues of $A$ and $B$. For the following proposition, denote by $\lambda_k(A)$, the $k$-th eigenvalue of $A$ with Dirichlet boundary condition when $k \geq 1$ and define $\lambda_0(A)=0$; use the same convention for $\lambda_k(B)$.
\begin{prop}\label{prop:Courant} With $A,B \subset M$ as above, let $k \geq 1$ and $0 \leq l \leq 2k$.  Then we have the inequality 
\begin{equation}\label{eq:nodal} \lambda_{2k-1}(M) \leq \max \big \{ \lambda_{2k-l}(A), \lambda_l(B) \big \}.\end{equation}
\end{prop}
For convenience and because we could not find a precise reference to this exact result in the literature, we give a short proof of Proposition \ref{prop:Courant} at the end of this subsection.

The Poincar\'e minimax principle states that 
	\begin{equation}\label{eq:Rayleigh}
		\lambda_k (D) = \inf_V \sup_{f \in V,\, f \not\equiv 0} \frac{\int_D \| \grad f \|^2 \, d\Vol_n}{\int_D f^2 \, d\Vol_n}
	\end{equation}
where $V$ runs over all
	\begin{enumerate}
		\item $k$-dimensional subspaces of $H_0^1(D)$, when $\partial D \neq \emptyset$,
		\item $(k+1)$-dimensional subspaces of $H^1(D)$, when $\partial D = \emptyset$.
	\end{enumerate}
\begin{rmk}
 The shift in dimension of $V$ is a consequence of the geometer's convention of indexing eigenvalues to start with $\lambda_0$ for closed manifolds.
\end{rmk}

For a discussion of the Sobolev spaces $H^1$ and $H^1_0$, see Appendix \ref{sec:Sobolev}.
 
\noindent {\bf Proof of Proposition \ref{prop:Courant}.} Define $V_A$ and $V_B$ as the following subspaces of $H_0^1(A)$ and $H_0^1(B)$ respectively:
\begin{align*}
	V_A:= &\{\text{span of the first }2k-l\text{ eigenfunctions on }A\} \subset H_0^1(A)\\
	V_B:= &\{\text{span of the first }l\text{ eigenfunctions on }B\} \subset H_0^1(B).
\end{align*}
Since functions in $V_A$ satisfy the Dirichlet boundary condition, functions in $V_A$ can be extended to functions in $H^1(M)$ by defining them to be zero on the complement of $A$. The analogous construction works for functions in $V_B$. These extensions allow us to construct $V=V_A \oplus V_B$ so that $V$ is a subspace of $H^1(M)$ with dimension $2k$.

In this proof, all integrals will be taken with respect to $d\Vol_n$.
Write
	\begin{equation*}
		R(f)= \frac{\int_M \|\grad f \|^2}{\int_M \|f\|^2}
	\end{equation*}
for the Rayleigh quotient on $M$.
 Then we have by the minimax principle that
\begin{align}\label{align:addfnc}
		\lambda_{2k-1} (M) \leq \sup_{f \in V} R(f) &= \sup_{f_A\in V_A, f_B \in V_B} \frac{\int_A \|\grad f_A\|^2 + \int_B \|\grad f_B\|^2}{\int_A \|f_A\|^2 + \int_B \|f_B\|^2} \\ \label{align:comb}
			&\leq \frac{\lambda_{2k-l}(A)\! \int_A \|f_A\|^2 + \lambda_l(B)\! \int_B \|f_B\|^2}{\int_A \|f_A\|^2 + \int_B \|f_B\|^2}\\
			&\leq \max \big \{\lambda_{2k-l}(A),\lambda_l(B) \big \}.
	\end{align}
The equality in (\ref{align:addfnc}) follows by writing $f=f_A+f_B$ where $f_A \in V_A$ and $f_B \in V_B$. Since $f_A$ is a linear combination of the first $2k-l$ eigenfunctions on $A$, its Rayleigh quotient over $A$ is at most $\lambda_{2k-l}(A)$. The analogous observation is also true for $f_B$, so its Raylaigh quotient at most $\lambda_l(B)$.  Therefore, the inequality (\ref{align:comb}) follows.
\eproof

\subsection{Single Parameter Test Functions on $M$}
 We now provide the setup for the proof of Theorem \ref{theo:SLE} giving upper bounds on $\lambda_k (M)$ in terms of an Sturm-Liouville problem which depends on $h(M)$. To do this, we first show how to give an upper bound for $\lambda_k (M)$ in terms of a Rayleigh quotient of a test function depending only on the distance to $\Sigma$. Our methods follow the arguments given in Buser \cite {B82} to obtain an upper bound for $\lambda_1(M)$ in terms of a Rayleigh quotient with a one-dimensional test function.

Recall that $\xi(k) := \left \lceil \frac{k+1}{2} \right \rceil$. Define 
	\begin{equation*}
		D_k:= \left \{ \begin{array}{lr} A, & \text{ if }\lambda_{\xi(k)} (A) \geq \lambda_{\xi(k)} (B),\\ B, & \text{ if }\lambda_{\xi(k)} (A) < \lambda_{\xi(k)} (B), \end{array} \right.
	\end{equation*}
for $k \in \N$.
Then by Proposition \ref{prop:Courant} with $l=k$ and Poincar\'e's minimax pinciple (\ref{eq:Rayleigh}), for a test function $f \in H_0^1(D_k)$, we have 
	\begin{equation}
		\lambda_k(M) \leq \lambda_{\xi(k)} (D_k)  = \inf_V \sup_{f \in V} \frac{\int_{D_k} \|\grad~f\|^2 \,d\Vol_n}{\int_{D_k} \|f\|^2\, d\Vol_n} \label{eq:Courant}
	\end{equation}
where $V$ ranges over $\xi(k)$-dimensional subspaces of $H_0^1(D_k)$.  
	\begin{rmk}To simplify notation, we will write $D_k$ as $D$, omitting the subscript $k$ where it is easily understood. In any case, the reader should remember that $D$ depends on $k \in \N$.\end{rmk}

Let $d:M \ra \R$ be the signed distance function given by 
$$d(x):=\left \{\begin{array}{rr}\dist(x, \Sigma=\partial D),& x \in D,\\ -\dist(x, \Sigma=\partial D),& x \notin D. \end{array}\right.$$
We now restrict the test functions $f$ in (\ref{eq:Courant}) to functions of the form $f=g\circ d$ where $g \in H_0^1[0,\infty)$. A posteriori, by Lemma \ref{lem:Main2}, it will be clear that we can take $g\in C_0^{\infty}[0, \infty)$. However, the following lemma shows that it is not necessary to make such a restriction.

\begin{lem}\label{lem:Chain}
	If $g\in H_0^1 [0,\infty)$ and $d:M \ra \R$ is the signed distance to $\Sigma$, then $g \circ d \in H^1(D)$.
\end{lem}

\begin{rmk}
	The standard chain rule for composition of Sobolev functions goes the other way around: the inner function is in $H^1$ and the outer function is Lipschitz, see for instance Evans and Gariepy \cite[Section 4.2, Theorem 4]{EG}. Example \ref{exam:tent} is a counter-example which shows that Lemma \ref{lem:Chain} is not true when $d$ is an arbitrary Lipschitz function.
\end{rmk}

\begin{exam}{\label{exam:tent}} Let $\Phi:[1,\infty) \ra \R$ be a smooth cut-off function with $\Phi(1)=1$ and $\Phi(x)=0$ when $x \geq 2$. Then let 
	\begin{equation*}
		g(x)= \left \{ \begin{array}[c]{rr} x^{\frac{3}{4}}, & 0\leq x \leq 1,\\ \Phi (x),& 1 < x <\infty, \end{array} \right.
	\end{equation*}
so $g \in H_0^1[0, \infty)$.

Define $d(x)$ as follows. Choose a sequence of numbers $a_j>0$ such that 
	\begin{equation}\label{eq:sequence}
		\sum_{j=1}^{\infty} a_j=1, \qquad \sum_{j=1}^{\infty} \sqrt{a_j}=+\infty.
	\end{equation}
For instance, we can let $a_j=cj^{-\frac{3}{2}}$ for suitable constant $c$. Let $\Lambda:\R \ra [0,1]$ be the following ``tent'' function of slope $\pm 1$, supported on $[-1,0]$:
	\begin{equation*}
		\Lambda (x)= \left \{ \begin{array}[c]{rr} x+1,& -1 \leq x \leq -\frac{1}{2},\\-x, & -\frac{1}{2} <x \leq 0,\\ 0,& 						\mathrm{~otherwise}. \end{array} \right.
	\end{equation*}
Let $s_j=a_1+\cdots +a_j$ be the $j$-th partial sum of the sequence, so $s_j \ra 1$, and define $d:\R \ra [0,\infty)$ such that 
	\begin{equation*}
		d(x)=\sum_{j=1}^{\infty} a_j \Lambda \! \left (\frac{x-s_j}{a_j} \right ).
	\end{equation*}
Then $d(x)$ is supported on $[0,1]$ and has slope $\pm 1$ at each point, except for the isolated local maximum and minimum points. Thus, $d$ is Lipschitz.

On the other hand, $g \circ d \notin H^1(\R)$, since
	\begin{equation*}
		(g \circ d)' = g'\big (d(x)\big ) \cdot d'(x) = \frac{3}{4}d(x)^{-\frac{1}{4}} \cdot (\pm 1).
	\end{equation*}
So then 
	\begin{align*}
		\int_{\R} (g \circ d)'(x)^2 \, dx &=\frac{9}{16} \int_0^1 d(x)^{-\frac{1}{2}} \, dx\\
						&=\frac{9}{16} \sum_{j=1}^{\infty} \int_0^1 a_j^{ -\frac{1}{2}} 
							\Lambda \left ( \frac{x-s_j}{a_j}\right )^{\!\!\!-\frac{1}{2}}\, dx\\
						&=\frac{9}{16} \sum_{j=1}^{\infty} \sqrt{a_j} 
						\left (\int_{-1}^0 \Lambda(y)^{-\frac{1}{2}} \, dy \right )\!\!, 													&& \text{by~the~change~of~variable}~
						y=\frac{x-s_j}{a_j},\\
						&=+\infty, &&\text{by~hypothesis}~(\ref{eq:sequence}).
	\end{align*}
\end{exam}

The proof of Lemma \ref{lem:Chain} will use the following Lemma.

\begin{lem}\label{lem:Chain2}
		Let $f \in L^2[0,\infty)$ and $d$ as above. Then there exists  a constant $C>0$ such that $\|f \circ d\|_{L^2(D)} \leq C \|f\|_{L^2[0,\infty)}$.
\end{lem}

\D {\bf Proof of Lemma \ref{lem:Chain2}.}
Since $d$ is a distance function to a rectifiable current and $M$ is compact, $\Vol_{n-1}\big (d^{-1}(\tau)\big )$ is bounded. Because $\|\grad(d)\|=1$ almost everywhere on $M$, the coarea formula gives
	\begin{align*}
		\int_D (f \circ d)^2 \, d\Vol_n &=\int_0^{\infty} \int_{d^{-1}(\tau)} (f \circ d)^2 \, d\Vol_{n-1} \, d\tau\\
			&=\int_0^{\infty}f^2(\tau)\Vol_{n-1}\big (d^{-1}(\tau)\big )\, d\tau\\
			&\leq C\|f\|^2_{L^2(D)}
	\end{align*}
where $C$ is an upper bound on the $(n-1)$-volume of the sets $d^{-1}(\tau)$.
\eproof

We now prove Lemma \ref{lem:Chain}.

\D {\bf Proof of Lemma \ref{lem:Chain}.} Suppose that $g\in C^{\infty}_c[0,\infty)$ and $X$ is a smooth vector field on $M$. Then by Rademacher's Theorem, since $d$ is Lipschitz, the derivative $X(d)$ exists almost everywhere. So then, by integration by parts, for all $\phi \in C^{\infty}_c(D)$,
	\begin{equation}\label{eq:ChainIBP}
		\int_D (g\circ d) X(\phi) \, d\Vol_n = - \int_D (g' \circ d) X(d) \phi \, d\Vol_n.
	\end{equation}
Thus, $g \circ d$ has a weak derivative $(g' \circ d) X(d)$. We will show that this weak derivative is square integrable.

Since $\|\grad(d)\|=1$ almost everywhere in $D$, the coarea formula gives
\begin{align*}
	\int_{D} (g' \circ d)^2 \, d\vol_n =& \int_0^{\infty} \int_{d^{-1}(\tau)} (g'\circ d)^2(\tau) \, d\Vol_{n-1}(\tau) \, d\tau\\
		=& \int_0^{\infty} (g')^2(\tau) \Vol_{n-1}\big (d^{-1}(\tau) \big ) \, d\tau< \infty.
\end{align*}
The last inequality follows from the the facts that $g'$ is a compactly supported function in $L^2[0,\infty)$ and $\Vol_{n-1}\big (d^{-1}(\tau)\big )$ is bounded and finitely supported on $[0,\infty)$. So $(g' \circ d) \in L^2(D)$. Further, 
	\begin{equation*}
		\|X(d)_p\| \leq \Lip(d)\|X_p\| \leq \|X_p\|
	\end{equation*}
for all $p \in M$ where the right hand side is uniformly bounded since $M$ is compact. Thus, we have \linebreak $(g' \circ d) X(d) \in L^2(D)$, and hence $g\circ d \in H_0^1(D)$.

Now consider an arbitrary $g \in H_0^1[0,\infty)=W^{1,2}_0[0,\infty)$ and approximate $g$ by a sequence of functions $g_k \in C_c^{\infty}(D)$ in the $W^{1,2}(D)$-norm. Then 
	\begin{equation*}
		\int_D \left ( |g-g_k|^2 + |g'-g_k'|^2 \right ) \, d\Vol_n \, \longrightarrow \, 0
	\end{equation*}
as $k \ra \infty$. 
So $g_k \circ d \ra g\circ d $ in $L^2(D)$ by Lemma \ref{lem:Chain2}, and $(g'_k \circ d) \ra (g' \circ d)$ in $L^2(D)$ by Lemma \ref{lem:Chain2}. Hence (\ref{eq:ChainIBP}) holds for $g$, by applying the result for $g_k\in C_c^{\infty}(D)$ and letting $k \ra \infty$. Thus, $g \circ d$ is weakly differentiable, with weak derivative in $L^2(D)$.
\eproof

We now resume bounding $\lambda_k(M)$ from above by a Rayleigh quotient with test functions whose values depend only on the distance to $\Sigma$. A routine calculation in Fermi coordinates shows that equation (\ref{eq:Courant}) implies 
	\begin{equation}
		\lambda_k(M) \leq \inf_{V_{\ast}} \sup_{g \in V_{\ast}} 
			\frac{\int_0^{\infty} \big (g'(\tau)\big )^{\! 2} \Vol_{n-1} \big (d^{-1}(\tau)\big )d\tau}{\int_0^{\infty}g^2(\tau) \Vol_{n-1} \big (d^{-1}(\tau)\big ) d\tau},
			\label{eq:Rayleigh1}
	\end{equation}
where we take $f=g \circ d$ and $V_{\ast}$ ranges over $\xi(k)$-dimensional subspaces of $H_0^1[0, \infty)$.

\subsection{Mean Curvature Bounds} 
\begin{figure}[htb]
\labellist
\small\hair 2pt
 \pinlabel {$h$} [ ] at 370 13
 \pinlabel {$\lambda(h)$} [ ] at 10 226
 \pinlabel {$\Eta(h)=\frac{h}{n-1}$} [ ] at 267 106
 \pinlabel {$\Eta(h)=1+\frac{h}{n}$} [ ] at 160 115
\endlabellist
\centering
\includegraphics[scale=1.0]{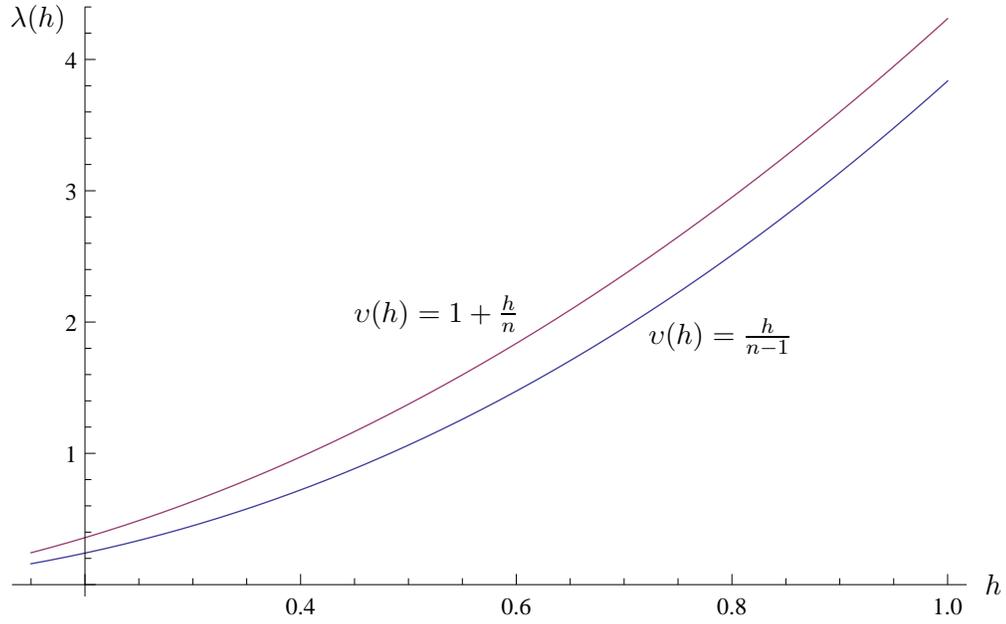}
\caption{Comparison of the upper bounds $\lambda(h)$ from Theorem \ref{theo:SLE} for $h \in [0.2, 1]$ for different choices of $\Eta(h)$.}
\label{fig:mean1}
\end{figure}

\begin{figure}[htb]
\labellist
\small\hair 2pt
 \pinlabel {$\lambda_1(h)$} [ ] at -5 220
 \pinlabel {$h$} [ ] at 370 13
 \pinlabel {$\Eta(h)=\frac{h}{n-1}$} [ ] at 210 87
 \pinlabel {$\Eta(h)=1+\frac{h}{n}$} [ ] at 145 110
\endlabellist
\centering
\includegraphics[scale=1.0]{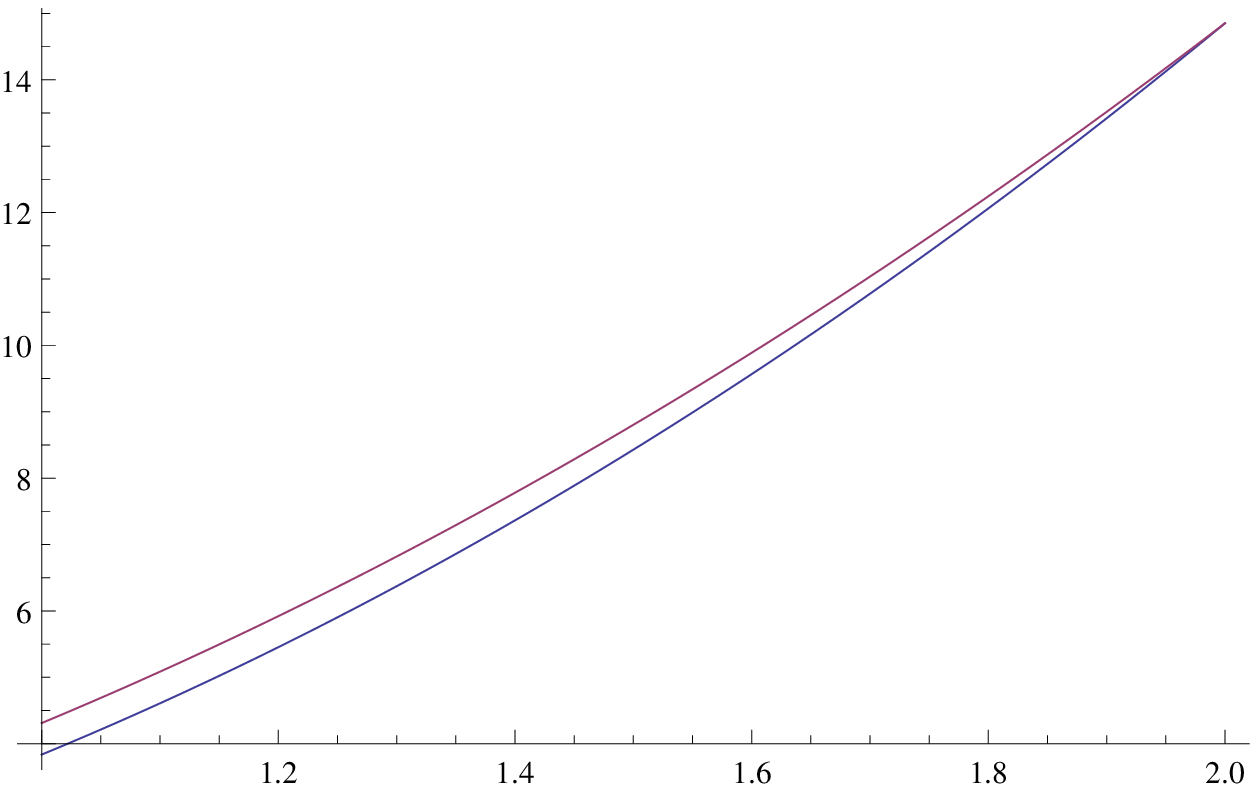}
\caption{Comparison of the upper bounds $\lambda_1(h)$ for $h \in [1, 2]$ for different choices of $\Eta(h)$.}
\label{fig:mean2}
\end{figure}

\begin{figure}[htb]
\labellist
\small\hair 2pt
 \pinlabel {$\lambda_1(h)$} [ ] at -5 227
 \pinlabel {$h$} [ ] at 370 13
 \pinlabel {$\Eta(h)=\frac{h}{n-1}$} [ ] at 178 123
 \pinlabel {$\Eta(h)=1+\frac{h}{n}$} [ ] at 270 113
\endlabellist
\centering
\includegraphics[scale=1.0]{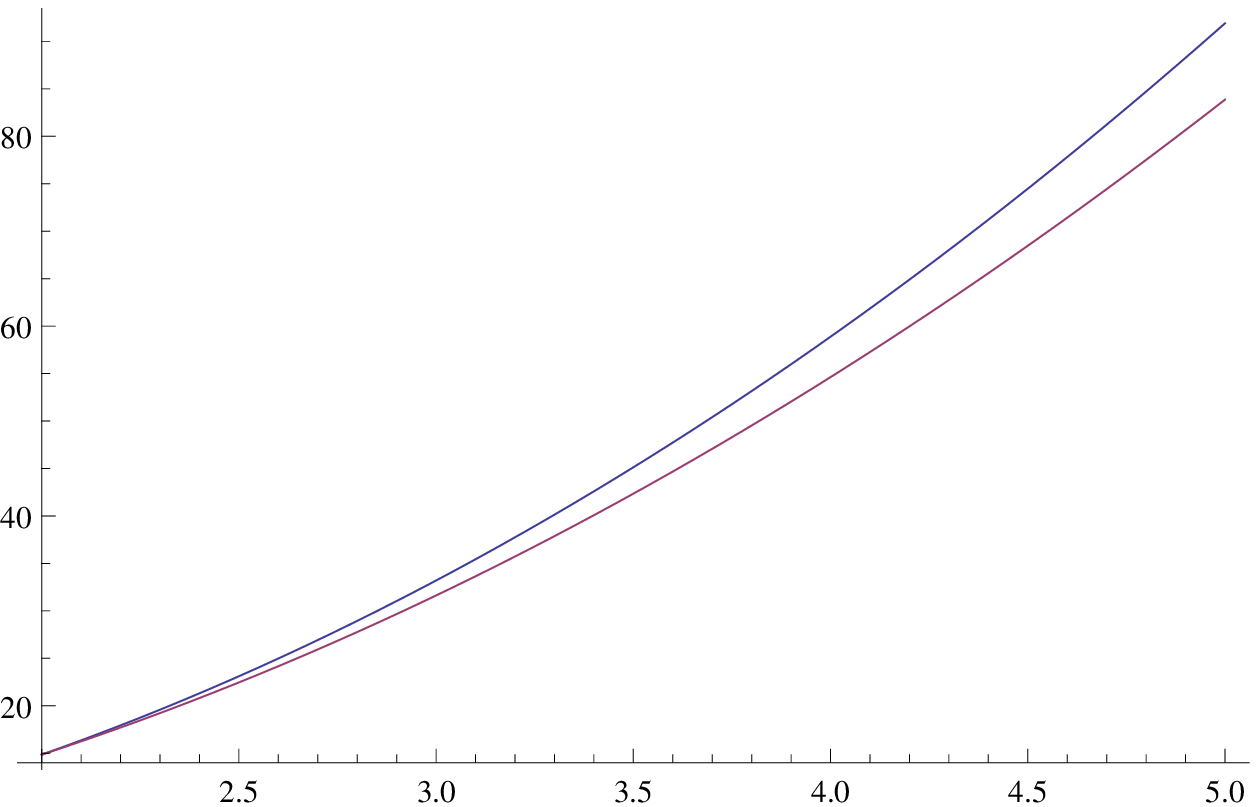}
\caption{Comparison of the upper bounds $\lambda_1(h)$ for $h \in [2, 5]$ for different choices of $\Eta(h)$.}
\label{fig:mean3}
\end{figure}

In order to further estimate the Rayleigh quotient for $\lambda_k(M)$, we consider a bound on the mean curvature of $\Sigma^0$, which is constant. Recall that this bound depends on $h(M)$ and we denote it by $\Eta (h)$. Buser's original approach used a comparison theorem of Heintze and Karcher \cite{HK78}, see Lemma \ref{lem:HK2}, to give an upper bound on the quantity $\Vol_{n-1}\big (d^{-1}(\tau)\big )$ in terms of the Ricci curvature and an upper bound on the mean curvature of $\Sigma$. Two simple upper bounds on the mean curvature of $\Sigma^0$ are $\Eta(h)=\frac{h}{n-1}$ given by Agol \cite{IA} and $\Eta (h)=\delta+\frac{h}{n}$ given by Buser \cite{B82}. Agol's bound has the benefit of not depending on the lower bound on Ricci curvature. In the case of $\delta=1$, Agol's bound is sharper when $h(M) <n(n-1)$ while Buser's bound on mean curvature is sharper when $h(M) >n(n-1)$. 
Figures \ref{fig:mean1}, \ref{fig:mean2}, and \ref{fig:mean3} give plots of the bounds $\lambda(h)$ for these two choices of $\Eta (h)$ for $n=2$ and $\delta=1$.

For $p \in \Sigma$, we denote $\eta (p)$ to be the mean curvature vector at $p$. The following statement was given by Agol \cite{IA}; we give a proof here for completeness.
\begin{lem} \label{lem:eta} If $\Sigma$ realizes the Cheeger constant and $\Vol_n (A) < \Vol_n (B)$, then $\eta$ points into $A$ everywhere.
\end{lem}

\noindent {\bf Proof of Lemma \ref{lem:eta}.}  First, proceed by contradiction assuming that $\Vol_n(A) < \Vol_n (B)$ and $\eta$ points into $B$.  Then there exists a current $\Sigma'$ which is a small perturbation of $\Sigma$ in the direction of $\eta$ at each point in $\Sigma$ with $\Sigma'$ separating $M$ into two disjoint regions $A'$ and $B'$ with $\Vol_n (A') < \Vol_n (B')$ with the convention that $A \subset A'$ and $B' \subset B$.  Since $\eta$ points into the direction of the perturbation, $\Vol_{n-1} (\Sigma') < \Vol_{n-1} (\Sigma)$.  Further, $\Vol_n (A) \leq \Vol_n (A')$, so then we have that
	\begin{equation*}
	h(M) = \frac{\Vol_{n-1}(\Sigma)}{\Vol_n (A)} > \frac{\Vol_{n-1} (\Sigma')}{\Vol_n (A')}.
	\end{equation*}
Since $\Vol_n(A') < \vol_n (B')$ implies that $\Vol_n (A') < \frac{1}{2} \Vol_n (M)$, we have a contradiction.
\eproof

The following result was given by Agol in order to give an upper bound for the norm of the mean curvature vector of $\Sigma^0$ in $M$. Since the mean curvature of $\Sigma^0$ is constant, we can refer to $H=\|\eta\|=\|\eta(p)\|$ for all $p \in \Sigma^0$ without ambiguity.
\begin{prop} {\bf (Agol \cite{IA})}  \label{prop:IA} For $\Sigma \subset M$ a Cheeger minimzing rectifiable current and $\Sigma^0 \subseteq \Sigma$ the smooth part of $\Sigma$, we have:
\begin{enumerate}
 	\item For $H$ on $\Sigma^0$, we have $(n-1) H \leq h(M)$.
	\item If $\Vol_n (A) < \Vol_n (B)$, we have $(n-1) H = h(M)$ and the mean curvature vectors $\eta(p)$ point into $A$ for all $p \in \Sigma^0$.
\end{enumerate}
\end{prop}

The following proof is the argument given by Agol in \cite{IA}.

\noindent {\bf Proof of Proposition \ref{prop:IA}.} Denote the cut locus of $\Sigma$ by $C=C(\Sigma)$.  Then Fermi coordinates on $M-C$ parameterize a subset $U \subset (-\infty, \infty) \times \Sigma$ with metric of the form $d\tau^2 + dA_{(\tau,s)}$ where $\tau$ is a signed distance from a point to $\Sigma$ and the $s$ is the minimizing geodesic point projection onto $\Sigma$. For more on the cut locus and Fermi coordinates, see Section \ref{sec:Fermi}.

Define $d_{\ast}:M \ra \R$ such that $$d_{\ast}(x):=\left \{\begin{array}{rr}\dist(x, \Sigma),& x \in A,\\ -\dist(x, \Sigma),& x \notin A. \end{array}\right. $$
Recall that $\Sigma^0$ is a hypersurface of constant mean curvature, so $H(p)$ is constant for all $p \in \Sigma^0$, and our convention is that $\Vol_n (A) \leq \Vol_n (B)$.  By Lemma \ref{lem:eta}, $\eta$ points into $A$. If $U \subset (-\infty, \infty) \times \Sigma$, then $d_{\ast}^{-1}(\tau)-C := \big ( \{ \tau \} \times \Sigma \big ) \cap U$.  Then the first variations of the volumes are
	\begin{equation}\label{eq:n1Vol}
		\left. \frac{d}{d\tau} \Vol_{n-1}\big ( d_{\ast}^{-1}(\tau) \big ) \right |_{\tau =0} = -(n-1) H \cdot \Vol_{n-1} (\Sigma)
	\end{equation}
and
	\begin{equation}\label{eq:nVol}
		\frac{d}{d\tau} \Vol_n \big (d_{\ast}^{-1}(\tau, \infty) \big ) = - \Vol_{n-1}\big (d_{\ast}^{-1}(\tau) \big ).
	\end{equation}
Applying the quotient rule and plugging in the first variation formulas (\ref{eq:n1Vol}) and (\ref{eq:nVol}), the infinitesimal change in the isoperimetric ratio is
	\begin{equation}\begin{split}\label{eq:Agol}
		\left. \frac{d}{d\tau} \frac{\Vol_{n-1}\big (d_{\ast}^{-1}(\tau)\big )}{\Vol_n \big (d_{\ast}^{-1}(\tau,\infty)\big )} \right |_{\tau=0} & =  
			\left ( \frac{\Vol_{n-1} \big (d_{\ast}^{-1}(\tau)\big )^2}{\Vol_n \big (d_{\ast}^{-1}(\tau, \infty)\big )^2} - \frac{(n-1)H 
				\Vol_{n-1} \big (d_{\ast}^{-1}(\tau)\big )}{\Vol_n \big (d_{\ast}^{-1}(\tau, \infty)\big )} \right )_{\!\!\! \tau=0}\\
			& =  h(h-(n-1) H).
	\end{split}\end{equation}

When $\Vol_n (A) < \Vol_n (B)$, we know that $\tau =0$ must be a critical point of $$\Vol_{n-1} \big (d_{\ast}^{-1}(\tau)\big )\left /\Vol_n \big (d_{\ast}^{-1}(\tau, \infty)\big )\right. .$$ If $\tau =0$ is not a critical point, we can perturb $\Sigma$ in the direction of $B$ which decreases $\Vol_{n-1} (\Sigma)$ and increases $\Vol_n (A)$, contradicting that $\Sigma$ realizes the Cheeger constant. So since $h \neq 0$, we have that $h=(n-1)H$.

When $\Vol_n (A) = \Vol_n (B)$, we can only consider the ratio 
\begin{equation*}
\Vol_n \left (d^{-1}_{\ast}(\tau)\right )\left /\Vol_n \left ( d^{-1}_{\ast}(\tau, \infty )\right )\right.
\end{equation*}
for $\tau>0$. If $\tau \leq 0$, we would have $\Vol_n(A) >\Vol_n(B)$ and the isoperimetric ratio $\frac{\Vol_{n-1}(\partial A)}{\Vol_n(A)}$ would not be a candidate for the Cheeger constant. So (\ref{eq:Agol}) gives us that $0>h\big (h-(n-1)H \big )$.  It follows that $h > (n-1) H$ in this case.
\eproof

We also consider the following mean curvature bounds given by Buser which depend on the lower bound of $-(n-1)\delta^2$ on the Ricci curvature of $M$.
\begin{prop}\label{prop:Buser} {\bf (Buser \cite{B82})}  With $M$, $\delta \geq 0$, and $H = |\eta |$ as above, then $H \leq \delta+ \frac{h}{n}$.
\end{prop}
The proof of Buser's result can be found at the end of the proof of Theorem 1.2 which can be found in Section 3 in \cite{B82}. We give a proof here for convenience. 

\noindent {\bf Proof of Proposition \ref{prop:Buser}.}  First we consider the case where $\delta >0$. Let $\tilde{J}(\tau, H) = \left (\cosh (\delta \tau) - \frac{H}{\delta}\sinh (\delta \tau) \right )^{n-1}$ when the term in parentheses is positive and $\tilde{J}(\tau, H) =0$ otherwise.  When $H \leq \delta$, the result follows immediately since $\frac{h}{n}\geq0$, so we may assume that $H > \delta$. Heintze and Karcher \cite{HK78} show that $$\int_0^{\infty}\tilde{J} (\tau , H)\, d\tau \geq 1/h.$$ We define $\epsilon>0$ so that $H = \delta(1+\epsilon)$ giving
	\begin{equation*} \cosh(\delta \tau ) - \frac{H}{\delta} \sinh(\delta \tau) = e^{- \delta \tau} - \epsilon \sinh (\delta \tau) \leq 1- \epsilon \delta \tau, \end{equation*}
noting that the right hand side is negative on $(1/\delta \epsilon, \infty)$.  It follows that
	\begin{align*} h= \frac{\Vol_{n-1} (\Sigma^0)}{\min \{\Vol_n (A), \Vol_n (B) \}} \geq& \frac{1}{\int_0^{\infty} \tilde{J}(\tau,H) \, d\tau}\\ \geq &\frac{1}{\int_0^{1/\delta \epsilon} (1-\epsilon\delta \tau)^{n-1}\, d\tau}\\ =& \epsilon \delta n = (H-\delta) n.
	\end{align*}
Letting $\delta \ra 0$ in the argument for $\delta >0$ gives the result for the case of $\delta=0$. 
\eproof

Combining Proposition \ref{prop:IA} of Agol and Proposition \ref{prop:Buser} of Buser, we arrive at the following observation.

\begin{prop}\label{prop:equal} Suppose that $M$ has Ricci curvature bounded below by $-(n-1)$. If $h(M) >n(n-1)$, then $\Vol_n (A)=\Vol_n(B)$.
\end{prop}

\D {\bf Proof of Proposition \ref{prop:equal}.}  Since $h>n(n-1)$, we have $1+\frac{h}{n} < \frac{h}{n-1}$.  Combining this observation with Buser's Proposition \ref{prop:Buser}, we have $H \leq 1 +\frac{h}{n} < \frac{h}{n-1}$.  If $\Vol_n (A) \neq \Vol_n (B)$, then $H=\frac{h}{n-1}$ by Agol's Proposition \ref{prop:IA}, a contradiction.  Thus, $\Vol_n (A) = \Vol_n (B)$.
\eproof

The next result follows from the work of Heintze and Karcher \cite{HK78} and was used as stated below by both Buser \cite{B82} and Agol \cite{IA}.

\begin{lem}\label{lem:HK2} {\bf (Heintze and Karcher \cite{HK78})} Let $d$ and $\Sigma$ be as previously defined and let $\Eta$ be a real number with $\Eta \geq H$. If the Ricci curvature of $M$ is bounded below by $-(n-1)\delta^2$, then 
\begin{equation}\label{eq:HK}
 \Vol_{n-1}\big (d^{-1}(\tau)\big ) \leq \Vol_{n-1} (\Sigma) J_{\delta}(\tau ,H) \leq \Vol_{n-1}(\Sigma)J_{\delta} (\tau, \Eta ).
\end{equation}
\end{lem}

We include a proof here for the convenience of the reader, following \cite{B82}.

\D {\bf Proof of Lemma \ref{lem:HK2}.}
Denote the solid tube of radius $R$ in the direction of the normal of $\Sigma^0$ by $\Tu (\Sigma, R)$.  Heintze and Karcher give
	\begin{equation}\label{eq:HK21}
		\Vol_n \big (\Tu (\Sigma, R) \big ) \leq \Vol_{n-1} (\Sigma) \int_{\xi \in S^0} \int_0^R \big ( c_{\delta} (\tau) - \langle \eta , \xi \rangle s_{\delta} (\tau) \big )^{n-1}\, d\tau \, d\xi
	\end{equation}
where the 0-sphere $S^0$ resides in $T_p\Sigma^{\perp}$ and the integrand is taken to be zero when $\big ( c_{\delta} (\tau) - \langle \eta , \xi \rangle s_{\delta} (\tau) \big )^{n-1}$ is negative \cite[Theorem 2.1]{HK78}.  Of the two vectors comprising $S^0$, one component points into $A$ and the other into $B$; denote these components $\xi^A_p$ and $\xi^B_p$ respectively and similarly for $\xi_p^D$.  Then the right hand side of (\ref{eq:HK21}) is equal to 
	\begin{equation*}
		\Vol_{n-1} (\Sigma) \left ( \int_0^R \big ( c_{\delta}(\tau) - \langle \eta, \xi^A_p \rangle s_{\delta} (\tau) \big )^{n-1} \, d\tau
			+ \int_0^R \big ( c_{\delta} (\tau) - \langle \eta, \xi^B_p \rangle s_{\delta} (\tau) \big )^{n-1} \, d\tau \right ).
	\end{equation*}
It follows that for $t \in (0,\infty)$,
	\begin{equation}\label{eq:ATUB1}
		\Vol_n\big (d^{-1}(0,t)\big ) \leq \Vol_{n-1} (\Sigma) \int_0^t \big ( c_{\delta} (\tau) - \langle \eta, \xi^D_p \rangle s_{\delta} (\tau)\big )^{n-1}\,  d\tau.
	\end{equation}
Now either $\langle \eta, \xi^D_p \rangle =H$ or $\langle \eta, \xi^D_p \rangle =-H$ for every $p \in \Sigma^0$.  Either way, from (\ref{eq:ATUB1}), we have
	\begin{equation}\label{eq:ATUB2}
		\Vol_n \big (d^{-1}(0,t)\big ) \leq \Vol_{n-1} (\Sigma) \int_0^t \big ( c_{\delta} (\tau) + H s_{\delta} (\tau) \big )^{n-1} \, d\tau.
	\end{equation}
Further, since $M$ is compact, the integrand on the right hand side is positive up until some value $t=t_1>0$, whereas for any $t_2 >t_1$, we have $d^{-1}(0,t_2)-d^{-1}(0,t_1)=\emptyset$. So $t>0$ can be as large as necessary and the inequality (\ref{eq:ATUB2}) will still hold.
This gives
	\begin{align*}
	 	\Vol_n \big (d^{-1}(0,t)\big )=\Vol_n \big (d^{-1}(0,t) - C\big ) &\leq  \Vol_{n-1} (\Sigma) \int_0^t J_{\delta} (\tau, H) d\tau\\ & \leq h(M) \Vol_n (D) \int_0^t J_{\delta} (\tau, H) d\tau
	\end{align*}
and $\Vol_{n-1} \big (d^{-1}(\tau)\big ) \leq \Vol_{n-1} (\Sigma ) J_{\delta}(\tau, H )$.
So then we have 
  \begin{equation*}
   J_{\delta} (\tau, H)= \big ( c_{\delta} (\tau ) + H s_{\delta} (\tau) \big )^{n-1} \leq \big ( c_{\delta} (\tau) + \Eta s_{\delta} (\tau) \big )^{n-1} = J_{\delta} (\tau, \Eta).
  \end{equation*}
This gives 
\begin{equation*}
 \Vol_{n-1}(d^{-1}(\tau)) \leq \Vol_{n-1} (\Sigma) J_{\delta} (\tau, H) \leq \Vol_{n-1}(\Sigma)J_{\delta}(\tau, \Eta).
\end{equation*}
\eproof

%%%%%%%%%%%%%%%%%%%%%%%%%%%%%%%%%%%%%%%%%%%%%%%%%%%%%%%%%%%%%%%%%%%%%%%%%%%%%%%%%%%%%%%%%%%%%%%%%%
\section{Distance Functions and Level Sets}
To prove the upper bounds on $\lambda_k(M)$ in terms of $h(M)$ for $k>1$ in Theorem \ref{theo:SLE}, we will need a lower bound on $\Vol_{n-1} \big (d^{-1}(\tau)\big )$ in terms of $\Vol_{n-1}(\Sigma)$. Recall the definition of $T$ from Section \ref{sec:detail} and that $D=A$ when $\lambda_k(A) \geq \lambda_k (B)$ and $D=B$ otherwise. We prove the following lemma:

\begin{lem}\label{lem:LSigTau} With $M, D, \Sigma,$ and $d$ as above,
	\begin{equation*}
	 \frac{\Vol_{n-1}\big (d^{-1}(\tau)\big )}{\Vol_n \big (d^{-1} (\tau, \infty)\big )} \geq \frac{\Vol_{n-1}(\Sigma)}{\Vol_n (D)}
	\end{equation*}
for $\tau \in (0,T)$ almost everywhere.
\end{lem}

We will prove Lemma \ref{lem:LSigTau} by proving three related lemmas. Specifically, Lemma \ref{lem:taugeqT} will help prove Lemma \ref{lem:taugeqT2}, while Lemmas \ref{lem:taugeqT2} and \ref{lem:cutlevel} will help prove Lemma \ref{lem:level}. Finally, Lemma \ref{lem:level} will be used to prove Lemma \ref{lem:LSigTau}.

\begin{lem}\label{lem:taugeqT}
If $\Vol_n \big (d^{-1} (\tau, \infty)\big )=0$ for $\tau >0$, then $\tau \geq T$.
\end{lem}

\noindent {\bf Proof of Lemma \ref{lem:taugeqT}.}
Suppose $\tau <T$. Then
	\begin{align*}
		\Vol_n\big (d^{-1}(0,\tau) \big )=\int_0^{\tau} \Vol_{n-1}\big (d^{-1}(x)\big )\, dx \leq & \int_0^{\tau} \Vol_{n-1} (\Sigma) J_{\delta}(x) \, dx\\ <& \int_0^T \Vol_{n-1}(\Sigma) J_{\delta}(x) \, dx \\= &\frac{\Vol_{n-1}(\Sigma)}{h}=\Vol_{n-1}(A) \leq \Vol_{n-1}(D),
	\end{align*}
by Lemma \ref{lem:HK2} and since $J_{\delta}(x)>0$ for $x\geq 0$. Thus, $\Vol_{n}\big (d^{-1}(\tau,\infty)\big )>0$.
\eproof

\begin{lem}\label{lem:taugeqT2}For any non-empty open set $E \subseteq (0,T)$, we have that $\int_E \Vol_{n-1} \big (d^{-1}(x)\big )\, dx >0.$ 
\end{lem}
\noindent {\bf Proof of Lemma \ref{lem:taugeqT2}.} By Lemma \ref{lem:taugeqT}, there must be a point of $D$ at least distance $T$ from $\Sigma$. Since $d$ is continuous, the interval $[0,T]$ is contained in $d(D)$; hence $d^{-1}(E)$ is a non-empty open subset of $D$. As such, it contains an open geodesic ball and, thus, $\int_E \Vol_{n-1} \big ( d^{-1}(x) \big )\, dx = \Vol_{n} \big ( d^{-1}(E) \big )>0$.
\eproof 

\begin{lem}\label{lem:cutlevel}
Let $C=C(\Sigma)$ be the cut locus of $\Sigma$ in $M$. Then $C \cap d^{-1}(\tau)$ has $(n-1)$-Hausdorff measure zero for $\tau \in \R$ almost everywhere.
\end{lem}

\D {\bf Proof of Lemma \ref{lem:cutlevel}.} Since $\| \grad(d) \|=1$ almost everywhere, the coarea formula gives:
	\begin{equation*}
		\Vol_n(C) = \int_{-\infty}^{\infty} \Vol_{n-1} \big ( C \cap d^{-1}(\tau) \big )\, d\tau.
	\end{equation*}
Therefore, since $\Vol_n(C)=0$, if follows that $\Vol_{n-1} \big (C \cap d^{-1}(\tau)\big )=0$ almost everywhere for $\tau \in \R$. \eproof

\begin{lem}\label{lem:level}
 Let $M$ and $T$ be as previously defined. Then for Lebesgue almost everywhere $\tau \in (0,T)$, $$h(M) \leq \frac{\Vol_{n-1}\big (d^{-1}(\tau)\big )}{\min \Big \{\Vol_n\big (d^{-1}(-\infty, \tau)\big ),\Vol_n\big (d^{-1} (\tau, \infty)\big )\Big \}}.$$
\end{lem}

\D {\bf Proof of Lemma \ref{lem:level}.} Because $\|\grad(d)\|=1$ almost everywhere on $M$, the slicing lemma tells us that $d^{-1}(\tau)$ is an $(n-1)$-rectifiable current for almost every $\tau \in (0,T)$; see Krantz and Parks \cite[Lemma 7.6.1]{KP} or Simon \cite[28.1 Lemma]{S83}. Since $d^{-1}(\tau)$ is the boundary of $d^{-1}[\tau, \infty)$, it follows that $d^{-1}[\tau,\infty)$ is an integral current for almost every $\tau \in (0,T)$. So then the Approximation Theorem, see Federer \cite[4.2.20]{HF69} and Morgan \cite[7.1]{FM09}, gives the following. For all $\epsilon >0$, there exists a finite simplicial complex $P$ which is smoothly embedded in $M$, such that $d^{-1}(\tau) = P+E$ where the current $E$ is such that $\Vol_{n-1}(E) < \epsilon$. It follows that $\Vol_{n-1}\big (d^{-1}(\tau)-P\big ) <\epsilon$. Further, since $P$ has codimension-1 in $M$, it is well-known that $P$ can be approximated by smooth submanifolds $S_{\delta_0}$ such that $\Vol_{n-1}\big (P -S_{\delta_0}\big )<\delta$. Then we 
have that $\Vol_{n-1}\big (\Sigma_{\tau}-S_{\delta_0}\big )<\epsilon + \delta_0$. Taking $\delta_0$ and $\epsilon$ to be arbitrarily small, by the definition of the Cheeger constant, we have that 
$$h(M) \leq \frac{\Vol_{n-1}\big (d^{-1}(\tau)\big )}{\min \Big \{\Vol_n\big (d^{-1}(-\infty, \tau)\big ),\Vol_n\big (d^{-1} (\tau, \infty)\big )\Big \}},$$
 since $\Vol_n\big (d^{-1}(-\infty, \tau)\big )$ and $\Vol_n\big (d^{-1} (\tau, \infty)\big )$ are strictly greater than or equal to zero for all $\tau \in (0,T)$ by Lemma \ref{lem:taugeqT2}.
\eproof

Now we will use Lemma \ref{lem:level} to prove Lemma \ref{lem:LSigTau}.

\D {\bf Proof of Lemma \ref{lem:LSigTau}.}
Here we apply Lemma \ref{lem:level}. Since $\Sigma$ has the property that $h(M)=\Vol_{n-1}(\Sigma)/\Vol_n(A)$, Lemma \ref{lem:level} gives that 
	\begin{equation}\label{eq:hbd}
	 \frac{\Vol_{n-1}(\Sigma)}{\Vol_n(A)} \leq \frac{\Vol_{n-1}\big (d^{-1}(\tau)\big )}{\min \Big \{ \Vol_n\big (d^{-1}(-\infty,\tau)\big ), \Vol_n \big (d^{-1}(\tau,\infty)\big ) \Big \}}
	\end{equation}
is true for Lebesegue $\tau \in (0,T)$ almost everywhere.
Working off of (\ref{eq:hbd}), we have two cases for almost every fixed $\tau \in (0,T)$.
  \begin{enumerate}
   \item First, we assume that $\Vol_n\big (d^{-1}(\tau, \infty)\big ) \leq \Vol_n\big (d^{-1}(-\infty,\tau)\big )$. Then we have
      \begin{equation*}
       \frac{\Vol_{n-1} (\Sigma)}{\Vol_n(D)} \leq \frac{\Vol_{n-1} (\Sigma)}{\Vol_n(A)} \leq \frac{\Vol_{n-1} \big (d^{-1}(\tau)\big )}{\Vol_n\big ( d^{-1}(\tau, \infty)\big )}.
      \end{equation*}
      So we have verified equation (\ref{eq:lowerbd}) in this case.
   \item We now assume that $\Vol_n\big (d^{-1}(\tau,\infty)\big ) > \Vol_n\big (d^{-1}(-\infty,\tau)\big )$. Then we have
      \begin{equation}\label{eq:Dcomp}
       \frac{\Vol_{n-1}(\Sigma)}{\Vol_n(A)} \leq \frac{\Vol_{n-1}\big (d^{-1}(\tau)\big )}{\Vol_n\big (d^{-1}(-\infty,\tau)\big )} = \frac{\Vol_{n-1}\big (d^{-1}(\tau)\big )}{\Vol_n\big (d^{-1}(\tau,\infty)\big )}.
      \end{equation}
But in this case, since $\Vol_n(D^C)+\Vol_n\big (d^{-1}(0,\tau)\big ) < \Vol_n(D)- \Vol_n\big (d^{-1}(0,\tau)\big )$, we have that $\Vol_n(D^C) < \Vol_n (D)$ because $\Vol_n\big (d^{-1}(0,\tau)\big )>0$. So we conclude that $D^C=A$ and $D=B$.
So $d$ takes $(0,\infty)$ to $B$. Then (\ref{eq:Dcomp}) becomes 
      \begin{equation*}
       \frac{\Vol_{n-1}(\Sigma)}{\Vol_n(A)} \leq \frac{\Vol_{n-1}\big (d^{-1}(\tau)\big )}{\Vol_n(A) + \Vol_n\big (d^{-1}(0,\tau)\big )}
      \end{equation*}
and we can multiply both sides by $\Vol_n(A)+\Vol_n\big (d^{-1}(0,\tau)\big )$ to obtain
      \begin{equation}\label{eq:SigTauLow}
        \left (1+ \frac{\Vol_n\big (d^{-1}(0,\tau)\big )}{\Vol_n(A)} \right )\Vol_{n-1} (\Sigma) \leq \Vol_{n-1}\big (d^{-1}(\tau)\big ).
      \end{equation}
It follows from (\ref{eq:SigTauLow}) that
	\begin{equation}\label{eq:STLow}
		\Vol_{n-1}(\Sigma) \leq \Vol_{n-1}\big (d^{-1} (\tau) \big ).
	\end{equation}
Combining (\ref{eq:STLow}) with the trivial fact that $\Vol_n(D) \geq \Vol_n \big (d^{-1}(\tau, \infty)\big )$, gives the result for this case.
  \end{enumerate}
\eproof

%%%%%%%%%%%%%%%%%%%%%%%%%%%%%%%%%%%%%%%%%%%%%%%%%%%%%%%%%%%%%%%%%%%%%%%%%%%%%%%%%%%%%%%%%%%%%%
\section{Upper Bounds as Eigenvalues of Sturm-Liouville Problems}
\label{sec:SLP}
In this section, we will give upper bounds for the spectrum of the Laplacian on $M$ in terms of the Cheeger constant using one-dimensional, regular Sturm-Liouville eigenvalue problems. In particular, the Rayleigh quotients for $\lambda_k(M)$ are bounded above by Rayleigh quotients of functions on certain compact intervals. The Rayleigh quotient for these functions uses weighted inner products where the weights depend on $h(M)$. We can then apply the spectral theorem to give the existence of the eigenvalues of each Rayleigh quotient and show that the corresponding eigenfunctions satisfy a regular Sturm-Liouville eigenvalue problem.

\subsection{Sturm-Liouville Problems}
We focus on Sturm-Liouville eigenvalue problems (or Sturm-Liouville problems) on the interval $(0,T)$.\footnote{We follow the convention of Zettl \cite{AZ}, writing Sturm-Liouville problems on the open interval $(0,T)$ even though the functions depend on the end points $0$ and $T$.} Our examples will consist of an operator of the form
\begin{equation*}%\label{eq:SLop}
			L_i=-\frac{1}{w_i(\tau)}\frac{d}{d\tau} \left ( p(\tau) \frac{d}{d\tau} \right ).
\end{equation*}
where $p$ and $w_i$ are the weight functions defined in Section \ref{sec:detail} as
	\begin{align*}
		p(\tau)&=w_1(\tau)=J_{\delta}(\tau)= \big ( c_{\delta}(\tau) + \Eta s_{\delta} (\tau) \big )^{n-1}, \\
		w_2(\tau)&=1 - h\int_0^{\tau} J_{\delta} \big (x, \Eta (h)\big ) \, dx.
	\end{align*}
Denote by $\omega_i(h)$ the Sturm-Liouville problem on $(0,T)$ given by 
	\begin{equation*}
		L_iu=\lambda u, \qquad u(0)=0, \qquad u'(T)=0.
	\end{equation*}

\subsection{Application of the Spectral Theorem}
In this section, we prove the following lemma which will help us prove Theorem \ref{theo:SLE}.

\begin{lem}\label{lem:Main2}
	For $i=1,2$, there exist eigenfunctions $\varphi_j$ which satisfy the Sturm-Liouville problem $\omega_i(h)$ so that $\varphi_j$ is smooth for each $j \in \N$ and $\varphi_j(T) \neq 0$. In addition,
		\begin{equation}\label{eq:SLRay}
			\lambda_j \big (\omega_i(h) \big ) = \frac{\int_0^T \left ( \varphi_j' \right )^2 p \, d\tau}{\int_0^T \left ( \varphi_j \right )^2 w \, d\tau}.
		\end{equation}
\end{lem}

To prove Lemma \ref{lem:Main2}, we will apply the version of the spectral theorem stated in Appendix \ref{sec:Spectral}. In doing so, we define the following Hilbert spaces which will correspond to our application of the spectral theorem. Let $\Hh_i=L^2([0,T], w_i \, d\tau)$ with inner product given by 
	\begin{equation*}
		(\psi_1, \psi_2)_{\Hh_i} = \int_0^{T} \psi_1 \psi_2 w_i \, d\tau.
	\end{equation*}
Further, define $\Kk_i=\{ \psi \in H^1([0,T]; p \, d\tau, w_i \, d\tau):\psi (0)=0\}$ and with inner product given by 
	\begin{equation*}(\psi_1, \psi_2)_{\Kk_i} = \int_0^{T} \big (\psi_1'\psi_2'p + \psi_1\psi_2 w_i\big )\, d\tau.
	\end{equation*}
Then $a_i(\psi_1, \psi_2)=(\psi_1, \psi_2)_{\Kk_i}$ is a bilinear, continuous, symmetric, and elliptic form from $\Kk_i \times \Kk_i$ to $\R$ for $i=1,2$.

\noindent {\bf Proof of Lemma \ref{lem:Main2}.}   We fix $i=1$ or $2$ and drop it from the notation, so e.g. $\Kk=\Kk_i$.  Note that $\Kk$ is continuously, densely, and compactly imbedded in $\Hh$. This follows from the classical imbedding of $H^1$ into $L^2$ and the equivalence of the $\Kk$-norm with the $H^1$-norm and the $\Hh$-norm with the $L^2$-norm since the weight functions $p$ and $w_i$ are positive almost everywhere on the compact interval $[0,T]$.
 Letting $\lambda_j := \gamma_j-1$ in the statement of the spectral theorem, given as Theorem \ref{theo:Spectral} in Appendix \ref{sec:Spectral}, gives the existence of an orthonormal basis $\{\varphi_j\}$ of weak eigenfunctions $\varphi_j \in \Kk$ satisfying
	\begin{equation}\label{eq:weakeigen}
		(\varphi_j,v )_{\Kk} = (\lambda_j + 1)(\varphi_j, v)_{\Hh}
	\end{equation}
for all $v \in \Kk$.

We now argue that the $\varphi_j$ satisfy the Sturm-Liouville equations and $\varphi_j(T) \neq 0$ and $\varphi_j'(T)=0$. Rewrite (\ref{eq:weakeigen}) as	
\begin{equation}\label{eq:weakeigen2}
	\int_0^T \varphi_j'v' p \, d\tau= \lambda_j \int_0^T \varphi_jvw \, d\tau.
\end{equation}
Because $p, w \in C^\infty[0,T]$, the elliptic regularity theorem guarantees that $\varphi_j \in C^{\infty}[0,T]$, so we can integrate the left side of (\ref{eq:weakeigen2}) by parts for all $v \in \Kk$. This gives
\begin{equation}\label{eq:weakBC} \left. \varphi_j'pv \right \vert_0^T - \int_0^{T} (\varphi_j'p)'v \, d\tau=\lambda_j \int_0^{T} \varphi_j w v \, d\tau. \end{equation}
Choosing $v$ to be in $H^1_0[0,T]$, 
	\begin{equation}\label{eq:weak1} \int_0^{T} - (\varphi_j'p)'v \, d\tau=\lambda_j \int_0^{T} \varphi_j w v \, d\tau \end{equation}
since we have $v(0)=v(T)=0$.
So then, (\ref{eq:weak1}) is equivalent to 
	\begin{equation}\label{eq:weak2} \int_0^{T} \left ((\varphi_j'p)'+\lambda_j \varphi_j w \right )v\, d\tau= 0.\end{equation}
Now (\ref{eq:weak2}) is true for all $v \in H_0^1[0,T]$ and $H_0^1[0,T]$ is dense in $L^2[0,T]$, so by approximation in $L^2[0,T]$, 
	\begin{equation} \label{eq:SLeig}L\varphi_j = \lambda_j \varphi_j\end{equation}
pointwise on $(0,T)$.

We now show that $\varphi_j$ satisfies the Neumann boundary condition at the right endpoint of $[0,T]$. We have just shown that the pointwise eigenvalue equation (\ref{eq:SLeig}) holds on $(0,T)$, so its weak form (\ref{eq:weakBC}) simplifies to show that $\left. \varphi_j'pv \right \vert_0^T=0$. Since $v(0)=0$, we have $\varphi_j'(T)p(T)v(T)=0$.
Choosing a $v \in \Kk$ with $v(T) \neq 0$ gives $\varphi_j'(T)p (T)=0$.  
Since $p (T) >0$ we must have the natural boundary condition $\varphi_j'(T)=0$.  

It remains to show that $\varphi_j(T) \neq 0$.  
Since $\varphi_j$ satisfies a second order linear ordinary differential equation with smooth coefficients, by existence and uniqueness of ordinary differential equations, if both $\varphi_j(T)=0$ and $\varphi_j'(T)=0$, then $\varphi_j \equiv 0$.  
But then $\varphi_j$ is not an eigenfunction, a contradiction.  
So we conclude that $\varphi_j (T) \neq 0$.  

The statement (\ref{eq:SLRay}) follows from combining (\ref{eq:weakeigen2}) with the following observations to conclude that $v=\varphi_j$: write $v=\sum_{j=1}^{\infty} (v,\varphi_j)_{\Hh} \varphi_j$ in (\ref{eq:weakeigen}), recall that the $\varphi_j$ are orthogonal to one another in $\Hh$, and then note the well-ordering of the $\lambda_j$ corresponding of $\varphi_j$. Since we have shown the equivalence of (\ref{eq:weakeigen2}) with the Sturm-Liouville problem $\omega (h)$, the result holds.
\eproof

\begin{rmk}
	Note that Theorem \ref{theo:SLE} holds when $T$ is replaced by any $t \in \R$ with $0<t<T$. Since any test function $g$ on $[0,t]$ can be extended to a test function on $[0,T]$ by $g(\tau)=g(t)$ for $\tau \in (t, T]$, one can conclude that
	 $$\frac{\int_0^T (g')^2p \, d\tau}{\int_0^T g^2w_i \, d\tau} \leq \frac{\int_0^t (g')^2p \, d\tau}{\int_0^t g^2w_i \, d\tau}.$$
\end{rmk}.

\subsection{Proof of Theorem \ref{theo:SLE}.} We begin with the case of $\lambda_1\big (\omega_1(h)\big )$.
We wish to minimize the Rayleigh quotient given in the expression (\ref{eq:Rayleigh1}) which is 
	\begin{equation}\label{eq:RaySig}
		\frac{\int_0^{\infty} \big (g'(\tau)\big )^2 \Vol_{n-1} \big (d^{-1}(\tau)\big ) \, d\tau}{\int_0^{\infty} g^2(\tau) \Vol_{n-1} \big (d^{-1}(\tau)\big ) \, d\tau}.
	\end{equation}
Restricting to test functions where $g (\tau ) =g(T)$ for all $\tau \geq T$, we have that (\ref{eq:RaySig}) is equal to 
	\begin{equation}\label{eq:RaySig2}
		\frac{\int_0^{T} \big (g'(\tau)\big )^2 \Vol_{n-1} \big (d^{-1}(\tau)\big ) \, d\tau}{\int_0^{T} \big (g^2 (\tau) - g^2(T)\big )\Vol_{n-1} (d^{-1}(\tau))\, d\tau + g^2(T) \Vol_{n} (D)}
	\end{equation}

Now we follow Buser in applying the Heintze-Karcher comparison theorem \cite{HK78}. In particular, we wish to compare equation (\ref{eq:RaySig2}) to the quotient
	\begin{equation}\label{eq:RaySig3}
		\frac{\int_0^{T} \big (g'(\tau)\big )^2 J_{\delta}(\tau, \Eta) \, d\tau}{\int_0^{T} g^2(\tau) J_{\delta}(\tau,\Eta) \, d\tau}.
	\end{equation}
By Heintze and Karcher \cite{HK78}, see our Lemma \ref{lem:HK2}, we have $$\Vol_{n-1} \big (d^{-1}(\tau)\big ) \leq \Vol_{n-1} (\Sigma ) J_{\delta}(\tau, \Eta ).$$

The eigenfunction $g=\varphi_1$ of the Sturm-Liouville problem $\omega_1(h)$ satisfies $(pg')'=-\lambda_1 wg$ on $(0,T)$ with $g(0)=0$. Theorem 0 in Everitt, Kwong, and Zettl \cite{EKZ} shows that since $J_{\delta}(\tau, \Eta)\geq 0$ for $\tau \in (0,T)$, the number of zeros of the eigenfunction corresponding to $\lambda_1$ of the quotient (\ref{eq:RaySig3}) is zero. Therefore, we may assume that $g \geq 0$ on $(0,T)$. Hence $(pg')' \leq 0$ and so $pg'$ is decreasing on $(0,T)$. Since $g'(T)=0$, we conclude that $pg' \geq 0$ on $(0,T)$, so $g' \geq 0$.

Because $g$ is monotone increasing on $(0,T)$, we have that $g^2(\tau) \leq g^2(T)$ for all $\tau \in [0,T]$, and so
	\begin{equation}\label{eq:denomg}
	\int_0^{T} \left (g^2(\tau)-g^2(T) \right ) \Vol_{n-1} (\Sigma_{\tau})\, d\tau \geq
		\Vol_{n-1}(\Sigma) \int_0^{T} \left (g^2(\tau) - g^2(T) \right ) J_{\delta}(\tau, \Eta) \, d\tau
	\end{equation}
by Lemma \ref{lem:HK2}. Thus, equation (\ref{eq:RaySig2}) is bounded above by 
	\begin{equation}\label{eq:RaySig4}
	 \frac{\Vol_{n-1}(\Sigma) \int_0^{T} \big (g'(\tau)\big )^2 J_{\delta}(\tau, \Eta) \, d\tau}{\Vol_{n-1}(\Sigma) \int_0^{T} \big (g^2 (\tau) - g^2(T)\big )J_{\delta}(\tau, \Eta) \, d\tau + g^2(T) \Vol_{n} (D)}.
	\end{equation}

Further, because $\Vol_n (A) \leq \Vol_n (D)$ and 
$$\frac{\Vol_n(A)}{\Vol_{n-1} (\Sigma)}=\frac{1}{h} = \int_0^{T} J_{\delta}(\tau, \Eta) \, d\tau,$$ 
we have that (\ref{eq:RaySig4}) is bounded above by
	\begin{equation}\label{eq:RaySig5}
		\frac{\int_0^{T} \big (g'(\tau)\big )^2 J_{\delta}(\tau, \Eta) \, d\tau}{\int_0^{T} \big (g^2 (\tau) - g^2(T)\big ) J_{\delta}(\tau, \Eta) \, d\tau+g^2(T) \frac{\Vol_{n} (A)}{\Vol_{n-1} (\Sigma)}} =
		\frac{\int_0^{T} \big (g'(\tau)\big )^2 J_{\delta}(\tau, \Eta) \, d\tau}{\int_0^{T} g^2(\tau) J_{\delta}(\tau,\Eta) \, d\tau}=\lambda_1 \big ( \omega_1(h) \big ).
	\end{equation}
  The result follows for $\lambda_1 (M)$ by the second statement in Lemma \ref{lem:Main2} since $J_{\delta}(\tau, \Eta)=p_1=w$.
  
For the case of $\lambda_k\big (\omega_1(h)\big )$, when $g$ does not correspond to the minimum non-zero value of the Rayleigh quotient (\ref{eq:RaySig}), we cannot guarantee that $g$ satisfies $g^2(\tau ) \leq g^2(T)$ for all $0 < \tau <T$ and hence (\ref{eq:denomg}) may not hold.\footnote{In fact, I computed many numerical examples of higher eigenfunctions which fail to have this property.} From Lemma \ref{lem:LSigTau}, we have that
	\begin{equation}\label{eq:lowerbd}
		\left ( 1-\frac{\Vol_{n}\big (d^{-1}(0,\tau)\big )}{\Vol_n (D)} \right ) \Vol_{n-1} (\Sigma) \leq \Vol_{n-1} \big (d^{-1}(\tau)\big )
	\end{equation}
for almost every $\tau \in (0, T)$. Further, by Lemma \ref{lem:HK2}, we have
	\begin{equation}\label{eq:Dtau}
		\Vol_n \big (d^{-1}(0,\tau)\big ) = \int_0^{\tau} \Vol_{n-1} \big (d^{-1}(x)\big ) \, dx \leq \Vol_{n-1}(\Sigma) \int_0^{\tau}  J_{\delta} (x, \Eta) \, dx.
	\end{equation}
It follows from (\ref{eq:Dtau}) that 
	\begin{equation}\begin{split}\label{eq:Dtau2}
		\frac{\Vol_{n} \big (d^{-1}(0,\tau)\big )}{\Vol_n (D)} &
		\leq \frac{\Vol_{n-1} (\Sigma) \int_0^{\tau} J_{\delta} (x,\Eta) \, dx}{\Vol_n (D)}\\
		 & \leq \frac{\Vol_{n-1} (\Sigma)}{\Vol_n (A)} \int_0^{\tau} J_{\delta} (x,\Eta) \, dx\\
		& = h \int_0^{\tau} J_{\delta} (x,\Eta) \, dx.
	\end{split}
	\end{equation}
Combining (\ref{eq:lowerbd}) and (\ref{eq:Dtau2}), we get 
	\begin{equation*}
		\left (1-h \int_0^{\tau} J_{\delta} (x, \Eta)\, dx \right )\Vol_{n-1} (\Sigma) \leq \Vol_{n-1}\big (d^{-1}(\tau)\big )
	\end{equation*}
almost everywhere on $(0,T)$. Therefore, we use this to decrease the denominator in equation (\ref{eq:RaySig}) to give the upper bound
	\begin{equation*}
		\frac{\int^{T}_0 \big (g'(\tau)\big )^2 J_{\delta} (\tau,\Eta) \, d\tau}{\int_0^{T} g^2(\tau) \left (1 - h\int_0^{\tau} J_{\delta} (x,\Eta) \, dx \right ) \, d\tau} = \frac{\int^{T}_0 \big (g'(\tau)\big )^2 p (\tau) \, d\tau}{\int_0^{T} g^2(\tau) w_2(\tau) \, d\tau}=\lambda_k \big ( \omega_2 (h) \big ).
	\end{equation*}
The result for $\omega_2(h)$ above holds by the second statement of Lemma \ref{lem:Main2}. This concludes the proof.
\eproof

\subsection{Upper Bounds of $\lambda_k$ as Functions of $h$} The Sturm-Liouville problem $\omega_i(h)$ is regular, meaning that $p,w_i>0$ with $p,w_i$ continuous functions on $[0,T]$. Further, $\omega_i(h)$ has homogeneous, self-adjoint separated boundary conditions, meaning that the boundary conditions can be written  as $N_1Y(0)+N_2Y(T)=0$ where $Y=\Big [\begin{smallmatrix} y\\ py' \end{smallmatrix} \Big ]$ and $N_1$ and $N_2$ are $2 \times 2$ matrices with real entries.
Since $w(x)L_i$ is a formally self-adjoint differential operator, $\omega_i(h)$ is separated and self-adjoint as an Sturm-Liouville problem. This classification of $\omega_i(h)$ allows us to apply several results from Sturm-Liouville theory to observe some properties of $\lambda_k\big (\omega_i(h)\big )$ when considered as a function of $h \in (0, \infty)$.

Recall that $p$, $w_i$, and $T$ depend on $h$, so we may write 
	$$\lambda_k\big (\omega_i (h)\big )=\lambda_k \Big (\omega_i \big (T(h), 1/p(h), w_i (h)\big )\Big )$$
where the left endpoint of the Sturm-Liouville problem is fixed at 0 and $T(h)$ is the right endpoint of the Sturm-Liouville problem.
The fact that $\lambda_k\big (\omega_i(h)\big )$ are real-valued follows from a result of Atkinson \cite{FA64} or from Everitt, Kwong, and Zettl \cite{EKZ}.
The conclusion that the $\lambda_k\big (\omega_i(h)\big )$ are continuous can be concluded from a result of Kong, Wu, and Zettl \cite{KWZ}, which gives the continuity of an Sturm-Liouville problem with respect to a single constraint, such as an endpoint of a weight function like $p$ or $w_i$, which is considered as a continuous variable of the Sturm-Liouville problem. Specifically, we have $1/p(\tau) \in L^1(0,T)$ and $w_i(\tau) \in L^1(0,T)$, so then the function $\lambda_k\big (\omega_i(T, 1/p, w)\big )$ is continuous in each component of $\omega_i$. Further $T, 1/p, w_i$ are continuous in the variable $h$, so it follows that $\lambda_k\big (\omega_i(h)\big )$ must be continuous in $h$.

To show the differentiability of $\lambda_k\big (\omega_i(h)\big )$ for $h \in (0,\infty)$ almost everywhere, one can consider results of Kong, Wu, and Zettl \cite{KWZ} and M\"oller and Zettl \cite{MZ96}. Specifically, by applying these results along with the chain rule to $\lambda_k\big (\omega_i(h)\big )$, one can give implicit formulas for $\frac{d\lambda_k}{dh}(\omega_i)$ in terms of these functions and normalized eigenfunctions of $\omega_i(h)$. We omit these details since we do not have a use for such a formula herein.

%%%%%%%%%%%%%%%%%%%%%%%%%%%%%%%%%%%%%%%%%%%%%%%%%%%%%%%%%%%%%%%%%%%%%%%%%%%%%%%%%%%%%%%%%%%%%%%
\section{The First Eigenfunction and the Riemann Differential Equation}
\label{sec:Riemann}

While there are many techniques for estimating eigenvalues of Sturm-Liouville problems, it is interesting to consider how the eigenfunctions for our Sturm-Liouville problems are related to other well-studied differential equations. This section is motivated by a comment of Agol in \cite{IA}, specifically that he did not know how to solve the differential equation in Theorem \ref{theo:Agol} except when $n=3$. We will prove Theorem \ref{theo:Main3}, which says that if $\varphi_1$ is the first non-zero eigenfunction with eigenvalue $\lambda=\lambda_1\big (\omega_1(h)\big )$, then the scaled function $y=\varphi_1p^{1/2}$ satisfies a Riemann differential equation which depends on $\lambda$. One can conclude that the branches of $y$ in $\C \cup \{\infty\}$ are given by hypergeometric functions. Further, the function $y$ herein is the same function which appears in Agol's Theorem \ref{theo:Agol} in \cite{IA}. To this end, we will provide a proof of Agol's Theorem \ref{theo:Agol}, although our proof appeals to 
Theorem \ref{theo:SLE}. 

\subsection{The Riemann Differential Equation} We now give some background on the Riemann differential equation which is related to Agol's differential equation in Theorem \ref{theo:Agol}. With respect to the notation, we will follow the conventions of Poole \cite{P60}. To define the Riemann differential equation, we consider distinct $a,b,c \in \C$ which will correspond to the regular singularities of the equation and we will denote their associated local exponents by $\{\alpha, \alpha'\}, \{\beta, \beta'\}, \{ \gamma, \gamma'\}$ where $\alpha+\alpha'+\beta+\beta'+\gamma+\gamma'=1$.
Define $q_1$ and $q_2$ by the following:
 	\begin{equation}\label{eq:q} 
		\begin{split}
			q_1(x) & :=  \left. (1-\alpha-\alpha')\right / (x-a)\\
					&\quad \left. +(1-\beta-\beta') \right / (x-b)\\
					&\quad \left. +(1-\gamma -\gamma') \right / (x-c),\\
		 	q_2(x) & :=  \big ( \alpha \alpha'(a-b)(a-c)(x-b)(x-c) \\
					&\quad + \beta \beta' (b-c)(b-a)(x-a)(x-c)\\ 
					& \quad \left.+ \gamma \gamma'(c-a)(c-b)(x-a)(x-b)\big ) \right /\\
					& \quad \quad (x-a)^2(x-b)^2(x-c)^2.
		\end{split}
	\end{equation}
Then the Riemann differential equation is given by 
\begin{equation}\label{eq:Riem}
	\frac{d^2y}{dx^2} +q_1(x)\frac{dy}{dx}+q_2(x)y=0.
\end{equation}
Solutions of (\ref{eq:Riem}) are branches of the corresponding Riemann $P$-function written as
	\begin{equation*}%\label{eq:sol}
		 P \left \{ \begin{array}[c]{cccc} a & c & b& \\ \alpha' & \gamma' & \beta' & x\\ \alpha & \gamma & \beta & \end{array} \right \}.
	\end{equation*} See Poole for additional details on the Riemann differential equation and the $P$-function \cite{P60}.

\subsection{Test Functions Satisfying a Riemann Differential Equation}
Let $q=q(n,\lambda)=\frac{n-1-2\lambda}{2}, r=r(n)=\frac{(n-1)(n-3)}{4},$ and $s=q+r$.  Further, let $a=\frac{\Eta +1}{\Eta -1}$, and $b=e^{2T}\frac{\Eta +1}{\Eta -1}$.  Recall the weight function introduced in Section \ref{section:intro}:
	 \begin{equation*} \label{eq:J} p(\tau, \Eta)=\big (\cosh(\tau)+\Eta (h)\sinh(\tau) \big )^{n-1},\end{equation*} for an upper bound on mean curvature $\Eta (h)$. Then we have that $\bar{J}=J^{1/2}=p^{1/2}$ for and $T \in (0, \infty)$ is defined implicitly by $$\frac{1}{h} = \int_0^{T} \bar{J}(\tau, \Eta) d \tau .$$ 

We will show that the scaled test function $y=g\bar{J}$ satisfies a Riemann differential equation on $[a,b]$. This allows us to give a proof for the reformulation of Agol's Theorem \ref{theo:Agol} given in Theorem \ref{theo:Main3}. Claims \ref{Claim:sing} and \ref{Claim:exp} and Lemma \ref{lem:Riem} below are all that are required to complete the proof of Theorem \ref{theo:Main3}.
\begin{Claim} \label{Claim:sing} The regular singularities of the differential equation (\ref{eq:Riem2}), $$y''(z) + \frac{1}{z} y'(z) - \frac{q(z-1)^2+r(z+1)^2}{4z^2(z-1)^2}y(z)=0$$ in $\C \cup \{\infty \}$, are $0,1, \infty$.  Further, (\ref{eq:Riem2}) is Fuchsian, i.e. all its singular points are regular singularities.
\end{Claim}

\D {\bf Proof of Claim \ref{Claim:sing}.}  Define $\rho_1(z), \rho_2(z)$ to be such that (\ref{eq:Riem2}) can be written in the form
	\begin{equation}\label{eq:p1p2}
	 	\frac{d^2y}{dz^2}+\rho_1(z)\frac{dy}{dz}+\rho_2(z)y=0.
	\end{equation}
It is well known that a point $P \in \C$ is not a regular singularity if and only if the $\rho_i$ in (\ref{eq:p1p2}) do not have a pole at $P$.  Further, a point $P \in \C$ is a regular singularity if and only if $\lim_{z \ra P}(z-P)^i \rho_i(z)$ exists for $i=1, 2$ and $\infty$ is a regular singularity if and only if $\lim_{z \ra \infty} z^i\rho_i(z)$ exists for $i=1,2$.  These details are given in references such as Beukers \cite{B07}.  Applying a partial fraction decomposition to $\rho_2(z)$ in (\ref{eq:p1p2}), we have
	\begin{equation}
		\rho_2(z)=-\left [\frac{r}{(z-1)^2} - \frac{r}{z-1}+\frac{s}{4z^2} + \frac{r}{z} \right ].\label{eq:p2}
	\end{equation}
Thus, from  (\ref{eq:p1p2}) and  (\ref{eq:p2}), the only poles of the $\rho_i$ are located at $0,1 \in \C$.  So it suffices to check the points $0,1, \infty \in \C \cup \{\infty\}$ for being regular singularities of (\ref{eq:ODE4}).  The limits $\lim_{z \ra P} (z-P)^i\rho_i(z)$ exist for $P=0,1$ and $i=1,2$, so $0$ and $1$ are regular singularities  of (\ref{eq:p1p2}).  Since $\rho_1(z)$ behaves like $\frac{1}{z}$ at $\infty$ and $\rho_2(z)$ behaves like $\frac{1}{z^2}$ at $\infty$, the limit of $z^i\rho_i(z)$ as $z$ goes to $\infty$ exists for $i=1,2$ and we conclude that $\infty$ is a regular singularity of (\ref{eq:p1p2}).\eproof

\begin{Claim} \label{Claim:exp} The local exponents of the regular singularities $0,1, \infty$ for equation (\ref{eq:Riem2}) are respectively $$\begin{array}{c}\left \{ -\frac{\sqrt{s}}{2}, \frac{\sqrt{s}}{2} \right \}, \left \{\frac{1-\sqrt{1+4r}}{2}, \frac{1+\sqrt{1+4r}}{2} \right \}, \left \{ -\frac{\sqrt{s}}{2}, \frac{\sqrt{s}}{2} \right \}.\end{array}$$
\end{Claim}

\D {\bf Proof of Claim \ref{Claim:exp}.} Recall that the local exponents of a regular singularity are the roots of the indicial equation corresponding to the singularity.
The form for the indicial equation for a singularity $t \in \C$ is given by
	\begin{equation*}%\label{eq:ind1}
		 X(X-1)+a_{1,t}X+a_{2,t}=0
	\end{equation*}
while the form of the indicial equation for the singularity at $\infty$ is
	\begin{equation*}%\label{eq:ind2} 
		X(X+1)-a_{1, \infty}X+a_{2, \infty}=0.
	\end{equation*}
Local exponents for each singularity are given by the roots of the respective indicial equation.
Using the following limits, the local exponents are given by routine calculations.  For $x$ at 0, we have
		\begin{equation*}
			 a_{1,0} = \lim_{x \ra 0} x\rho_1(x) =  1, ~
			 a_{2,0} = \lim_{x \ra 0} x^2\rho_2(x) =  -\frac{s}{4}
		\end{equation*}
For $x$ at 1, we have
		\begin{equation*}
			 a_{1,1} = \lim_{x \ra 1} (x-1)\rho_1(x) = 0, ~ 
			 a_{2,1} = \lim_{x \ra 1} (x-1)^2\rho_2(x) = -r
		\end{equation*}
For $x$ at $\infty$, we have
		\begin{equation*}
			 a_{1, \infty} = \lim_{x \ra \infty} x\rho_1(x) = 1, ~
			 a_{2, \infty} = \lim_{x \ra \infty} x^2\rho_2(x) =-\frac{s}{4}
		\end{equation*}
The claim follows. \eproof

\begin{lem} \label{lem:Riem} The equation $$y''(z) + \frac{1}{z} y'(z) - \frac{q(z-1)^2+r(z+1)^2}{4z^2(z-1)^2}y(z)=0,$$  can be realized as an example of Riemann's differential equation.  As a consequence, branches of the corresponding Riemann $P$-function solve the differential equation.
\end{lem}

\D {\bf Proof of Lemma \ref{lem:Riem}.} We must show that  $q_1(z)= \frac{1}{z}$ and $q_2(z)=  - \frac{q(z-1)^2+r(z+1)^2}{4z^2(z-1)^2}$.  Since $\infty$ is a singularity of (\ref{eq:Riem2}), we must take the limit as the corresponding singularity in the formulation of (\ref{eq:Riem}) goes to $\infty$.  Thus, with an abuse of notation, we let $c= \infty$.

Now since $c$ is linear in both the numerator and denominator of $q_1(x)$, we have that
	\begin{equation}\label{eq:limq1}
		\begin{split}
		 q_1(x) &=  \lim_{c \ra \infty} \frac{(1-\alpha-\alpha')(x-b)(-c)+(1-\beta-\beta')(x-a)(-c)}{(x-a)(x-b)(-c)}\\
		&= \frac{(1-\alpha-\alpha')(x-b)+(1-\beta-\beta')(x-a)}{(x-a)(x-b)}.
		\end{split}
	\end{equation}
Similarly, $c$ is quadratic in both the numerator and denominator of $q_2(x)$, the same argument as in (\ref{eq:limq1}) gives
	\begin{equation*} %\label{eq:limq2}
		 q_2(x) =  \frac{\alpha \alpha' (a-b)(x-b) + \beta \beta' (b-a)(x-a)
 		+ \gamma \gamma' (x-a)(x-b)}{(x-a)^2(x-b)^2}.
	\end{equation*}

Choose singularities $a=0$ and $b=1$ in (\ref{eq:q}).  Letting
	\begin{equation*}%\label{eq:abg}
		\begin{array}[c]{ll}		
			\alpha =-\frac{\sqrt{s}}{2}, & \alpha' =\frac{\sqrt{s}}{2},\\ 
			\beta =\frac{1-\sqrt{1+4r}}{2}, & \beta'= \frac{1+\sqrt{1+4r}}{2},\\
			\gamma = -\frac{\sqrt{s}}{2}, & \gamma' = \frac{\sqrt{s}}{2},
		\end{array}
	\end{equation*}
a routine calculation gives the result.
\eproof

We are now ready to prove Theorems \ref{theo:Agol} and \ref{theo:Main3}.

\subsection{Proof of Theorem \ref{theo:Agol} and Theorem \ref{theo:Main3}.} \label{sec:AgolRiem} We will show that these theorems follow from Theorem \ref{theo:SLE}; consider the differential equation $-(u'p)'=\lambda u w_1$ corresponding to $\omega_1(h)$. Since $\bar{J}=p^{1/2}$, then we can write $-(u'p)'=\lambda u w_1$ as $-(u'\bar{J}^2)'=\lambda u\bar{J}^2$. Letting $u=y/\bar{J}$, we have that $$y''=\left (\frac{\bar{J}'}{\bar{J}} - \lambda (h) \right )y.$$ The Dirichlet condition $y(a)=0$ is given by the Dirichlet boundary condition on $u$; in other words, $u(0)=0$.  The normalization $y(b)=\bar{J}(T)$ follows from the fact that since $u(T) \neq 0$ in the Sturm-Liouville equation, we can normalize $u$ so that $u(T)=1$.  Finally, the Neumann condition $y'(b)=\bar{J}'(T)/2b$ follows from the Neumann condition $u'(T)=0$ in the Sturm-Liouville equation since
	\begin{equation*}
		u'(T) = \frac{\bar{J}(T)y'(T)- y(T) \bar{J}'(T)}{\bar{J}^2(T)} =0.
	\end{equation*}
This completes a proof of Agol's Theorem \ref{theo:Agol}.

Now suppose that our mean curvature bound on $\Sigma$ satisfies $\Eta \neq 1$.  It is a routine computation to see that the Sturm-Liouville equation gives
	\begin{equation*}%\label{eq:ODE2}
		 y'' - \left [q+r \left ( \frac{(\Eta +1)e^{2\tau}+(\Eta -1)}{(\Eta +1)e^{2\tau}-(\Eta -1)} \right )^{\! 2} \right ] y= 0. 
	\end{equation*}
Note that $(\Eta +1)e^{2\tau} - (\Eta -1) \neq 0$ since $\tau \geq 0$.
Making the substitution $z=e^{2\tau}\frac{\Eta +1}{\Eta -1}$, where $\Eta \neq 1$, a routine computation gives
	\begin{equation}\label{eq:ODE4} 
		\frac{d^2y}{dz^2}+ \rho_1(z)\frac{dy}{dz}+\rho_2(z)y=0.
	\end{equation}
where $\rho_1(z)=\frac{1}{z}$ and $\rho_2(z) = - \frac{q(z-1)^2+r(z+1)^2}{4z^2(z-1)^2}$.

Recall that $s=q+r$.  It follows from Claims \ref{Claim:sing} and \ref{Claim:exp} and Lemma \ref{lem:Riem} that the differential equation (\ref{eq:ODE4}) is of the form of Riemann's differential equation, and thus, the solutions of (\ref{eq:ODE4}) are given by branches of the Riemann $P$-function. This completes the proof. \eproof

%%%%%%%%%%%%%%%%%%%%%%%%%%%%%%%%%%%%%%%%%%%%%%%%%%%%%%%%%%%
\section{Examples and Applications}
\label{sec:exam}

\subsection{Examples} In practice, one might wish to estimate $\lambda_k(M)$ from above using $h(M)$. We consider three simple examples where both quantities are known explicitly and consider an Sturm-Liouville problem $\omega (h)$ given by finding an appropriate scaling function $J$ and give $\lambda_k\big (\omega (h)\big )$.
\begin{exam}
The 1-sphere, $S^1$, with the metric induced by embedding it in $\R^2$ with radius 1. The eigenvalues are of the form $\lambda(S^1)=j^{2}$ with multiplicity $j+1$ for $j \in \N$. 

Here $\Sigma$ is the set of two antipodal points on $S^1$. Further, $\Vol_0(\Sigma)=\Vol_0\big (d^{-1}(\tau)\big )=2$ for $\tau \in \left [ 0, \frac{\pi r}{2}\right )$, so we can take $J(\tau)=1$ to be our scaling function for $\Vol_0 \big (d^{-1}(\tau) \big )$ in terms of $\Vol_0(\Sigma)$. Now since $h(S^1)=\frac{2}{\pi}$, we have $T=\frac{1}{h}=\frac{\pi}{2}$. This gives the following Sturm-Liouville problem $\omega$ on $\left ( 0, \frac{\pi}{2}\right )$:
	\begin{equation*}
		-(u')' = \lambda u, \quad u(0)=0,\quad u' \left (\frac{\pi}{2} \right )  =0.
	\end{equation*}
Since solutions of $\omega$ are of the form $u(\tau)=C\sin\big ( (1+2j) \tau \big )$ for $j \in \N$ and $C \in \R$, we have $\lambda_k (\omega)=\left ( 1+2k \right )^{2}$. While $\lambda_1(M)=1$ and $\lambda_1\big (\omega (h)\big )=9$, the higher eigenvalues both grow like $k^2$.
\end{exam}

\begin{exam}
The 2-sphere, $S^2$, with the metric induced by embedding it in $\R^3$ with radius 1. Eigenvalues are of the form $\lambda(S^2) = j(j+1)$ with multiplicity ${j+2 \choose 2}-{j \choose 2}$ for $j \in \N$. 

Here $\Sigma$ is a great circle and $\Vol_1(\Sigma)=2\pi$. Further, $\Vol_1\big (d^{-1}(\tau)\big ) = 2\pi \cos \left ( \tau \right )$ for $\tau \in \left [ 0, \frac{\pi}{2}\right )$. Taking $J(\tau)=\cos \left ( \tau \right )$, we have that $\Vol_1\big (d^{-1}(\tau)\big )=J(\tau)\Vol_1(\Sigma)$. Since $h(S^2)=1$, it follows that $T=\sin^{-1}\! \left (\frac{1}{h}\right )=\frac{\pi }{2}$. This gives the following Sturm-Liouville problem $\omega$ on $\left (0, \frac{\pi}{2} \right )$:
	\begin{equation*}
		-\big (\cos (\tau ) u' \big )' = \lambda \cos ( \tau ) u,\qquad
		u(0)=0,\qquad
		u' \left (\frac{\pi}{2} \right )  =0.
	\end{equation*}
Using SLEIGN2, the eigenvalues of this Sturm-Liouville problem are $$\lambda=2,12,30,56,90,132,182,240,306,380, 462, 552, 650,\ldots$$
\end{exam}

\begin{exam}
The $n$-torus $\torus^n = \R^n /\Z^n$. Let $\vec{v} \in \Z^n$, then eigenvalues are of the form $\lambda(\torus^n)=\big (2\pi \|\vec{v}\|\big )^2$ where each eigenvalue has multiplicity $n$.

Here $\Sigma$ can be given by the planes $\big \{(x_1, \ldots, x_n):x_n=0\big \}$ and $\left \{(x_1, \ldots, x_n):x_n=\frac{1}{2}\right \}$ in the fundamental region $[0,1) \times \cdots \times [0,1)$ in $\R^n$. So then we have $h(\torus^n)=4$ and $\Vol_{n-1}(\Sigma)=\Vol_{n-1}\big (d^{-1}(\tau)\big )=2$ for each distance $\tau \in \left [ -\frac{1}{4}, \frac{1}{4} \right )$ off of $\Sigma$, so we take $J(\tau)=1$. It follows that $T=\frac{1}{h}=\frac{1}{4}$. This gives the following Sturm-Liouville problem $\omega$ on $\left (0, \frac{1}{4} \right )$:
	\begin{equation*}
		-(u')' = \lambda u,\qquad
		u(0) =0,\qquad
		u' \left (\frac{1}{4} \right ) =0.
	\end{equation*}
Since solutions of $\omega$ are of the form $u(\tau)=C\sin \big (2\pi(1+2j)\tau\big )$ for $j \in \N$ and $C \in \R$, we have $\lambda_k(\omega)=\big (2\pi(1+2k)\big )^2$.
\end{exam}

\subsection{When $M$ is not Symmetric About $\Sigma$}
In examples where the geometry of $M$ is not symmetric about $\Sigma$, it may be possible to take a scaling function $J$ which is also not symmetric about $\Sigma$. Using the methods of Sections \ref{sec:eigen} and \ref{sec:SLP}, such a $J$ can be used to give two distinct Sturm-Liouville problems, $\omega_A$ so that $\lambda_k(A) \leq \lambda_k(\omega_A)$ and $\omega_B$ so that $\lambda_k(B) \leq \lambda_k(\omega_B)$.

Denote by $\spec(A)$, the set of eigenvalues with multiplicites of the Laplacian on $A$ with Dirichlet BC on $\partial A=\Sigma$ and let $\spec(B)$ be defined similarly for $B$. Then, it is a simple consequence of the Poincar\'e principle that 
	\begin{equation} \label{eq:SpecU}
		\lambda_{2k-1}(M) \leq \inf_S \sup_{\lambda \in S} ~\spec(A) \bigsqcup \spec(B)
	\end{equation}
where $S \subset  \spec(A) \bigsqcup \spec(B)$ with $|S|=2k$. In examples where $\spec (A) \neq \spec (B)$, it is straight-forward to see that applying (\ref{eq:SpecU}) in place of (\ref{eq:nodal}) in Proposition \ref{prop:Courant}, for some choice of $l$, gives sharper upper bounds for some of the values $\lambda_k(M)$.

%%%%%%%%%%%%%%%%%%%%%%%%%%%%%%%%%%%%%%%%%%%%%%%%%%%%%%%%%%%%%%%%%%%%%%
\section{Appendix}\label{section:appen}
This section contains a brief description of background material used herein.

\subsection{Sobolev Spaces}\label{sec:Sobolev}  We remind the reader of basic definitions of Sobolev spaces.  For multi-index $\alpha$ and function $u$, let $D_{\alpha}u$ denote the weak derivative of $u$ with respect to $\alpha$.
 Recall that the Sobolev space $H^k(M)$ is the completion of the set $\{u \in C^{\infty}(M): \|u\|_{k,2} < \infty\}$ where
	$$\|u\|_{k,2} = \sum_{0\leq |\alpha |\leq k} \left ( \int_M |D_{\alpha} u|^2\, dx \right )^{\!\!\frac{1}{2}}.$$
Further, $W^{k,2}(M)$ are functions $u \in L^2(M)$ such that $D_{\alpha}u$ exists and belongs to $L^2(M)$ for all $0\leq |\alpha | \leq k$.  Meyers and Serrin \cite{MS64} showed that $H^k(M)=W^{k,2}(M)$ and so we may use the two descriptions interchangeably as is convenient.  For a manifold with boundary $D$, we take $H^k_0(D)$ to be the functions $f \in H^k(D)$ such that $\left. f \right |_{\partial D} \equiv 0$, or equivalently, as the completion of the set $C^{\infty}_c(M)$ with respect to the Sobolev norm. For additional background on Sobolev spaces, see Evans and Gariepy \cite{EG} or Hebey \cite{EH99}.

\subsection{Background on Rectifiable Currents}\label{sec:current} An $(n-1)$-dimensional rectifiable current $\Sigma$ is an oriented subset of $M$ that is rectifiable in the Hausdorff measure $\Ha_{n-1}$.  That is, $\Sigma \subset M$ is a countable union of Lipschitz images of bounded subsets $B^{\ast}$ in $\R^{n-1}$ with $\Ha_{n-1} (B^{\ast}) < \infty$, ignoring sets of Hausdorff $(n-1)$-measure 0.  When we say $\Sigma$ has compact support, we think of $\Sigma$ as a function on $(n-1)$-forms $\phi$ given by
	$$\Sigma (\phi) = \int_{\Sigma} \langle \vec{S} (x), \phi (x) \rangle \mu (x) d \Ha_{n-1}$$
where $\vec{S}$ is the unit normal vector associated with the oriented tangent plane to $\Sigma$ at $x$ and $\mu (x)$ is an integer multiplicity, a nonnegative, integer-valued function with 
	$$\int_{\Sigma} \mu (x) d\Ha_{n-1} < \infty.$$  For the currents we study, $\mu(x)= 1$ for all $x \in \Sigma$, so that the mass of $\Sigma$ is exactly the $(n-1)$-volume of $\Sigma$. Further, $\langle \vec{S} (x), \phi (x) \rangle=\phi(x) \left ( \vec{S} (x)\right )$ is the pairing of the differential form $\phi(x)$ applied to the vector $\vec{S} (x)$.
The reader might wish to consult Federer \cite{HF69} and Morgan \cite{FM09} or Simon \cite{S83} for an overview of rectifiable currents.

\subsection{Fermi Coordinates and the Cut Locus} \label{sec:Fermi} We remind the reader of the definition of the cut locus of $\Sigma$ and the associated Fermi coordinates. The cut locus $C=C(\Sigma)$ of $M$ with respect to $\Sigma$ is the closure of the set of points $q \in M$, such that either
	\begin{enumerate}
		\item there exist two or more distance minimizing geodesics from $q$ to $\Sigma$ or
		\item $q$ is conjugate to a point in $\Sigma$ along a geodesic which joins them.
	\end{enumerate}
  It is well-known that $C$ is a closed set of Lebesgue measure zero, see for instance Gallot, Hulin, and Lafontaine \cite{GHL} or Cheeger \cite{C90}.  Now points $x \in M-C$ lie on a unique distance minimizing geodesic $\gamma:[0,1] \ra M$ from a point $x_{\Sigma} \in \Sigma$.  This geodesic points in the direction of the normal vectors to $\Sigma$ at the point $x_{\Sigma}$, so long as the normal vectors to $\Sigma$ are nonzero within the local neighborhood.  The point $x$ can then be represented by the parameters $(x_{\Sigma}, \tau)$, called Fermi coordinates, where $\tau$ is the distance between $x$ and $x_{\Sigma}$ along the geodesic $\gamma$.

\subsection{Hilbert Spaces and the Spectral Theorem}\label{sec:Spectral}
To prove the upper bounds on $\lambda_k (M)$ in terms of the eigenvalues of Sturm-Liouville problems, we use the spectral theorem. Recall that if $\Hh$ and $\Kk$ are infinite-dimensional Hilbert spaces over $\R$, then $\Kk$ is continuously and densely imbedded in $\Hh$ if there exists a continuous linear injection $\iota : \Kk \ra \Hh$ with $\iota (\Kk)$ dense in $\Hh$.  Further, an imbedding is compact if every bounded sequence in $\Kk$ has a subsequence that converges in $\Hh$.  Finally, a bilinear form $b: \Kk \times \Kk \ra \R$ is called elliptic if there exists a $C >0$ such that $b( \psi, \psi ) \geq C \| \psi \|_{\Kk}^2$ for all $\psi \in\Kk$.  For additional information on the spectral theorem, the reader might wish to consult Chapter 6 in Blanchard and Br\"uning \cite{BB92}. The following statement of the spectral theorem is used to prove our generalization of Agol's theorem.

\begin{theo}\label{theo:Spectral} {\bf (Spectral Theorem)} Let $\Hh$ and $\Kk$ be infinite dimensional Hilbert spaces over $\R$ with $\Kk$ continuously and densely imbedded in $\Hh$ and let this imbedding be compact.  Let $(\cdot , \cdot )_{\Kk}: \Kk \times \Kk \ra \R$ be bilinear, continuous, symmetric, and elliptic. 

Given these assumptions, there exist vectors $\varphi_1, \varphi_2, \varphi_3, \ldots \in \Kk$ and numbers $0 < \gamma_1 \leq \gamma_2 \leq \gamma_3 \leq \cdots \ra \infty$ such that 
\begin{itemize}
	 \item $\varphi_j$ is an eigenvalue of $(\cdot , \cdot )_{\Kk}$ with eigenvalue $\gamma_j$, i.e. for all $v \in \Kk$, 
			\begin{equation}
				a( \varphi_j , v)= \gamma_j (\varphi_j, v)_{\Hh},
			\end{equation}
	\item $\{ \varphi_j \}$ is an orthonormal basis for $\Hh$, and 
	\item $\{\varphi_j / \sqrt{\gamma_j} \}$ is an orthonormal basis for $\Kk$.
\end{itemize}
Finally, the decomposition
	\begin{equation}
		f= \sum_j (f,\varphi_j)_{\Hh} \cdot \varphi_j
	\end{equation}
converges in $\Hh$ for all $f \in \Hh$, and in $\Kk$ for all $f \in \Kk$.
\end{theo}

%\backmatter

\bibliographystyle{plain}
\bibliography{SLCheegSpec}

\begin{thebibliography}{10}

\bibitem{IA}
Ian Agol.
\newblock An improvement to {B}user's inequality.

\bibitem{FA76}
F.~J. Almgren, Jr.
\newblock Existence and regularity almost everywhere of solutions to elliptic
  variational problems with constraints.
\newblock {\em Mem. Amer. Math. Soc.}, 4(165):viii+199, 1976.

\bibitem{FA64}
F.~V. Atkinson.
\newblock {\em Discrete and continuous boundary problems}.
\newblock Mathematics in Science and Engineering, Vol. 8. Academic Press, New
  York-London, 1964.

\bibitem{AM87}
F.~V. Atkinson and A.~B. Mingarelli.
\newblock Asymptotics of the number of zeros and of the eigenvalues of general
  weighted {S}turm-{L}iouville problems.
\newblock {\em J. Reine Angew. Math.}, 375/376:380--393, 1987.

\bibitem{BEZ}
P.B. Bailey, W.N. Everitt, and A.~Zettl.
\newblock The {SLEIGN}2 {S}turm-{L}iouville code.
\newblock {\em ACM Trans. Math. Software}, 21:143--192, 2001.

\bibitem{BBG85}
P.~B{\'e}rard, G.~Besson, and S.~Gallot.
\newblock Sur une in\'egalit\'e isop\'erim\'etrique qui g\'en\'eralise celle de
  {P}aul {L}\'evy-{G}romov.
\newblock {\em Invent. Math.}, 80(2):295--308, 1985.

\bibitem{B07}
Frits Beukers.
\newblock Gauss' hypergeometric function.
\newblock In {\em Arithmetic and geometry around hypergeometric functions},
  volume 260 of {\em Progr. Math.}, pages 23--42. Birkh\"auser, Basel, 2007.

\bibitem{BB92}
Philippe Blanchard and Erwin Br{\"u}ning.
\newblock {\em Variational methods in mathematical physics}.
\newblock Texts and Monographs in Physics. Springer-Verlag, Berlin, 1992.
\newblock A unified approach, Translated from the German by Gillian M. Hayes.

\bibitem{B82}
Peter Buser.
\newblock A note on the isoperimetric constant.
\newblock {\em Ann. Sci. \'Ecole Norm. Sup. (4)}, 15(2):213--230, 1982.

\bibitem{C69}
Jeff Cheeger.
\newblock A lower bound for the smallest eigenvalue of the {L}aplacian.
\newblock In {\em Problems in analysis ({P}apers dedicated to {S}alomon
  {B}ochner, 1969)}, pages 195--199. Princeton Univ. Press, Princeton, N. J.,
  1970.

\bibitem{C90}
Jeff Cheeger.
\newblock Critical points of distance functions and applications to geometry.
\newblock In {\em Geometric topology: recent developments ({M}ontecatini
  {T}erme, 1990)}, volume 1504 of {\em Lecture Notes in Math.}, pages 1--38.
  Springer, Berlin, 1991.

\bibitem{C75a}
Shiu~Yuen Cheng.
\newblock Eigenfunctions and eigenvalues of {L}aplacian.
\newblock In {\em Differential geometry ({P}roc. {S}ympos. {P}ure {M}ath.,
  {V}ol. {XXVII}, {S}tanford {U}niv., {S}tanford, {C}alif., 1973), {P}art 2},
  pages 185--193. Amer. Math. Soc., Providence, R.I., 1975.

\bibitem{C75b}
Shiu~Yuen Cheng.
\newblock Eigenvalue comparison theorems and its geometric applications.
\newblock {\em Math. Z.}, 143(3):289--297, 1975.

\bibitem{EG}
Lawrence~C. Evans and Ronald~F. Gariepy.
\newblock {\em Measure theory and fine properties of functions}.
\newblock Studies in Advanced Mathematics. CRC Press, Boca Raton, FL, 1992.

\bibitem{EKZ}
W.~N. Everitt, Man~Kam Kwong, and A.~Zettl.
\newblock Oscillation of eigenfunctions of weighted regular {S}turm-{L}iouville
  problems.
\newblock {\em J. London Math. Soc. (2)}, 27(1):106--120, 1983.

\bibitem{HF69}
Herbert Federer.
\newblock {\em Geometric measure theory}.
\newblock Die Grundlehren der mathematischen Wissenschaften, Band 153.
  Springer-Verlag New York Inc., New York, 1969.

\bibitem{HF70}
Herbert Federer.
\newblock The singular sets of area minimizing rectifiable currents with
  codimension one and of area minimizing flat chains modulo two with arbitrary
  codimension.
\newblock {\em Bull. Amer. Math. Soc.}, 76:767--771, 1970.

\bibitem{GHL}
Sylvestre Gallot, Dominique Hulin, and Jacques Lafontaine.
\newblock {\em Riemannian geometry}.
\newblock Universitext. Springer-Verlag, Berlin, third edition, 2004.

\bibitem{Grom2}
Misha Gromov.
\newblock {\em Metric structures for {R}iemannian and non-{R}iemannian spaces},
  volume 152 of {\em Progress in Mathematics}.
\newblock Birkh\"auser Boston, Inc., Boston, MA, 1999.
\newblock Based on the 1981 French original [ MR0682063 (85e:53051)], With
  appendices by M. Katz, P. Pansu and S. Semmes, Translated from the French by
  Sean Michael Bates.

\bibitem{EH99}
Emmanuel Hebey.
\newblock {\em Nonlinear analysis on manifolds: {S}obolev spaces and
  inequalities}, volume~5 of {\em Courant Lecture Notes in Mathematics}.
\newblock New York University, Courant Institute of Mathematical Sciences, New
  York; American Mathematical Society, Providence, RI, 1999.

\bibitem{HK78}
Ernst Heintze and Hermann Karcher.
\newblock A general comparison theorem with applications to volume estimates
  for submanifolds.
\newblock {\em Ann. Sci. \'Ecole Norm. Sup. (4)}, 11(4):451--470, 1978.

\bibitem{KWZ}
Qingkai Kong, Hongyou Wu, and Anton Zettl.
\newblock Dependence of the {$n$}th {S}turm-{L}iouville eigenvalue on the
  problem.
\newblock {\em J. Differential Equations}, 156(2):328--354, 1999.

\bibitem{KP}
Steven~G. Krantz and Harold~R. Parks.
\newblock {\em Geometric integration theory}.
\newblock Cornerstones. Birkh\"auser Boston, Inc., Boston, MA, 2008.

\bibitem{K92}
Pawel Kr{\"o}ger.
\newblock On the spectral gap for compact manifolds.
\newblock {\em J. Differential Geom.}, 36(2):315--330, 1992.

\bibitem{K97}
Pawel Kr{\"o}ger.
\newblock On explicit bounds for the spectral gap on compact manifolds.
\newblock {\em Soochow J. Math.}, 23(3):339--344, 1997.

\bibitem{L51}
J.P. L{\'e}vy.
\newblock Probl{\'e}mes concrets d'nalyse fonctionnelle.
\newblock Paris, 1951.

\bibitem{LY80}
Peter Li and Shing~Tung Yau.
\newblock Estimates of eigenvalues of a compact {R}iemannian manifold.
\newblock In {\em Geometry of the {L}aplace operator ({P}roc. {S}ympos. {P}ure
  {M}ath., {U}niv. {H}awaii, {H}onolulu, {H}awaii, 1979)}, Proc. Sympos. Pure
  Math., XXXVI, pages 205--239. Amer. Math. Soc., Providence, R.I., 1980.

\bibitem{MS64}
Norman~G. Meyers and James Serrin.
\newblock {$H=W$}.
\newblock {\em Proc. Nat. Acad. Sci. U.S.A.}, 51:1055--1056, 1964.

\bibitem{LM}
Laurent Miclo.
\newblock On hyperboundedness and spectrum of markov operators.
\newblock \url{http://hal.archives-ouvertes.fr/ hal-00777146/}, 2013.
\newblock [HAL : hal-00777146, version 3].

\bibitem{MZ96}
Manfred M{\"o}ller and Anton Zettl.
\newblock Differentiable dependence of eigenvalues of operators in {B}anach
  spaces.
\newblock {\em J. Operator Theory}, 36(2):335--355, 1996.

\bibitem{FM03}
Frank Morgan.
\newblock Regularity of isoperimetric hypersurfaces in {R}iemannian manifolds.
\newblock {\em Trans. Amer. Math. Soc.}, 355(12):5041--5052 (electronic), 2003.

\bibitem{FM09}
Frank Morgan.
\newblock {\em Geometric measure theory}.
\newblock Elsevier/Academic Press, Amsterdam, fourth edition, 2009.
\newblock A beginner's guide.

\bibitem{P60}
E.~G.~C. Poole.
\newblock {\em Introduction to the theory of linear differential equations}.
\newblock Dover Publications, Inc., New York, 1960.

\bibitem{R05}
Antonio Ros.
\newblock The isoperimetric problem.
\newblock In {\em Global theory of minimal surfaces}, volume~2 of {\em Clay
  Math. Proc.}, pages 175--209. Amer. Math. Soc., Providence, RI, 2005.

\bibitem{S83}
Leon Simon.
\newblock {\em Lectures on geometric measure theory}, volume~3 of {\em
  Proceedings of the Centre for Mathematical Analysis, Australian National
  University}.
\newblock Australian National University, Centre for Mathematical Analysis,
  Canberra, 1983.

\bibitem{AZ}
Anton Zettl.
\newblock {\em Sturm-{L}iouville theory}, volume 121 of {\em Mathematical
  Surveys and Monographs}.
\newblock American Mathematical Society, Providence, RI, 2005.

\end{thebibliography}

\end{document}